\documentclass[11pt]{article}
\usepackage{amssymb,amsmath,amsthm,latexsym,amscd,amsfonts,amscd,graphics}
\usepackage[dvips]{graphicx}

\usepackage{latexsym}
\usepackage{makeidx}
\usepackage{multicol}
\usepackage{color}

\newtheorem{proposition}{Proposition}[]

\newtheorem{definition}{Definition}[]
\newtheorem{remark}{Remark}[]

\title{Diamond Twin}

\author{Toshikazu Sunada}

\date{}
\begin{document}

\maketitle

\begin{abstract}
As noticed in 2006 by the author of the present article, the hypothetical crystal---described by crystallographer F. Laves (1932) for the first time and designated ``Laves' graph of girth ten" by geometer H. S. M. Coxeter (1955)---is a unique crystal net sharing a remarkable symmetric property with the diamond crystal, thus deserving to be called the {\it diamond twin} although their shapes look quite a bit different at  first sight. In this short note, we shall provide an interesting mutual relationship between them, expressed in terms of ``building blocks" and ``period lattices." This may give further justification to employ the word ``twin." What is more, our discussion brings us to the notion of ``orthogonally symmetric lattice," a generalization of {\it irreducible root lattices}, which makes the diamond and its twin very distinct among all crystal structures. 
\end{abstract}

\begin{figure}[htbp]
\vspace{-0.5cm}
\begin{center}
\hspace{1cm}
\includegraphics[width=.62\linewidth]{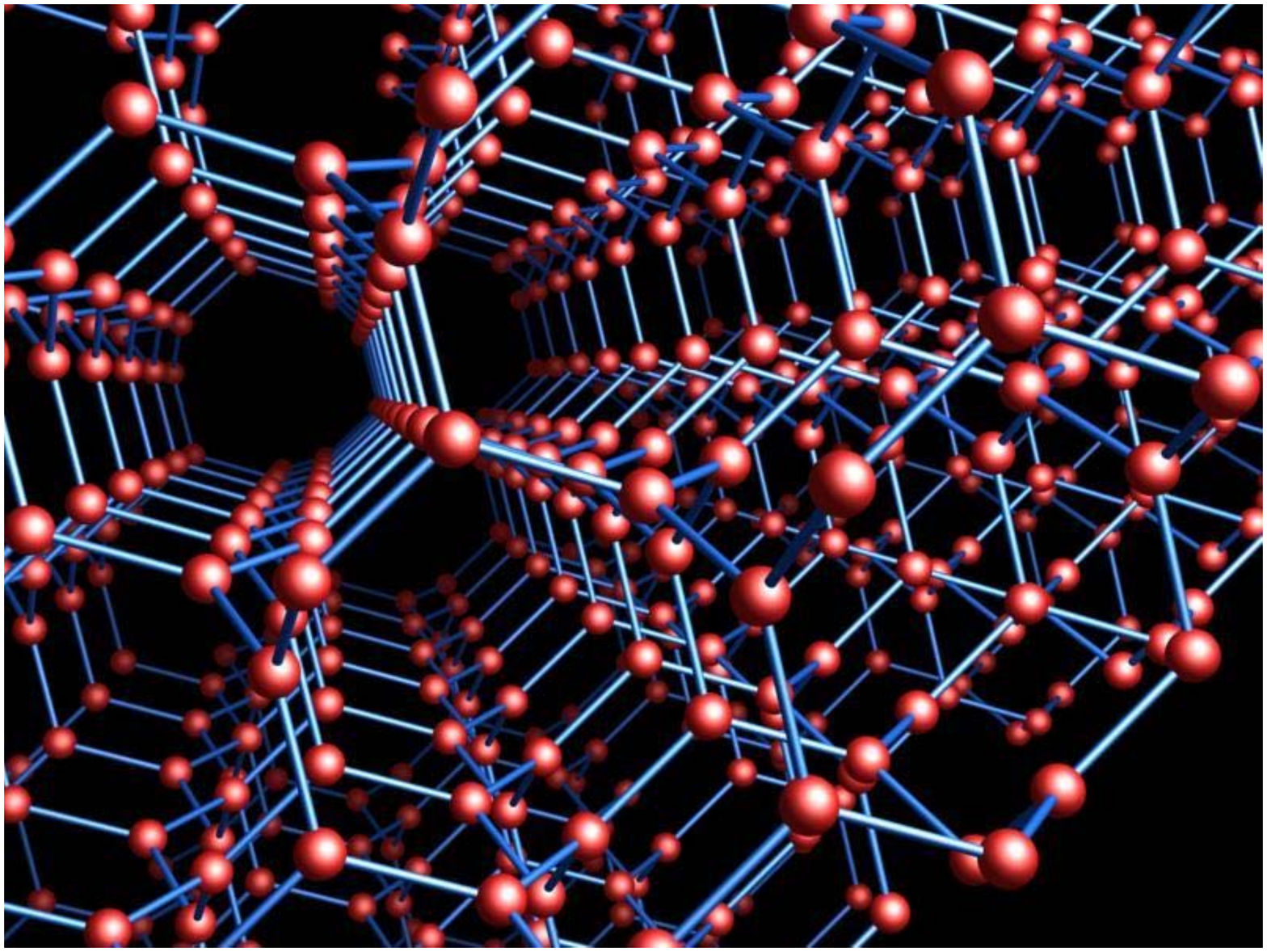}
\end{center}
\vspace{-0.6cm}
\caption{Diamond twin (CG-image created by Kayo Sunada)}\label{fig:k4Kayo}
\end{figure}

\section{Introduction}

The network structure of the diamond twin (Fig.~\!\ref{fig:k4Kayo}), which was described by German crystallographer Fritz Laves \cite{laves} for the first time as a hypothetical crystal and rediscovered by chance in the midst of my study of random walks on periodic graphs (\cite{kot}, \cite{su2}, \cite{su1}), is mathematically defined as the {\it standard realization} of the maximal abelian covering graph of the complete graph $K_4$ with 4 vertices (see \cite{kotani2}, \cite{kotani1}, \cite{su2} for the terminology).\footnote{We shall disregard the size of crystal nets and the physical character of atoms and atomic force since we are chiefly concerned with geometric properties. See \cite{su5} for the possible physical aspect of the diamond twin composed of carbon atoms. See also \cite{dai}, \cite{lian}, \cite{liu}, \cite{mizu} for related matters.} Actually this structure has been discussed and rediscovered by many people, and hence bears several names; {\it Laves' graph of girth ten} (\cite{cox}), ${\bf (10,3)}$-${\bf a}$ (\cite{wells}), the {\bf srs} net (the name coming from the fact that it occurs in the compound ${\rm SrSi}_2$ as the Si substructure; \cite{dd1}), the $K_4$ {\it crystal} for the obvious reason (\cite{su2}, \cite{su1}), and the {\it triamond net} (\cite{conway}). In addition, it has a close relationship with the {\it gyroid}, an infinitely connected triply periodic minimal surface discovered by Alan Schoen in 1970 (\cite{shoen}).\footnote{See \cite{hyde} for an interesting story concerning Laves' graph.}

\begin{figure}[htbp]
\vspace{-0.7cm}
\hspace{-0.7cm}
\includegraphics[width=1.08\linewidth]{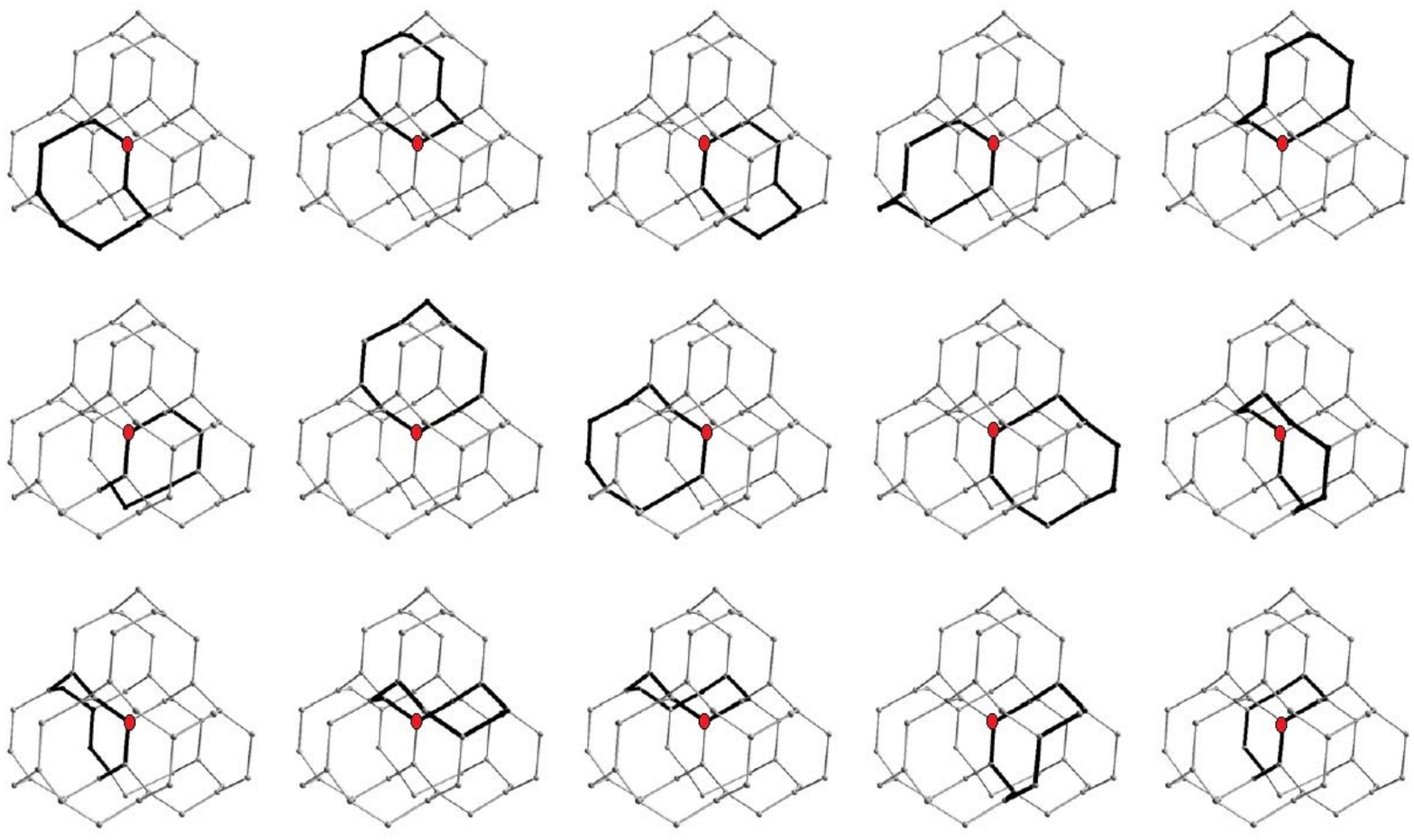}
\vspace{-2.0cm}
\caption{15 decagonal rings (CG-image created by Hisashi Naito)}\label{15}
\end{figure}

The reason why Laves' graph attracts much attention is that it possesses many remarkable properties from a geometrical point of view:

\begin{enumerate}

\item It has {\it maximal symmetry} in the sense that every automorphism of the diamond twin as an abstract graph extends to a congruent transformation of space (note that any congruent transformation fixing a crystal net induces an automorphism).

\item It has the {\it strongly isotropic property};  meaning that for any two vertices $x$ and $y$ of the crystal net, and for any ordering of the directed edges with the origin $x$ and any ordering of the directed edges with the origin $y$, there is a net-preserving congruence taking $x$ to $y$ and each $x$-edge to the similarly ordered $y$-edge (\cite{su1}).\footnote{We should not confuse the strongly isotropic property with the {\it edge-transitivity} or the notion of {\it symmetric graphs}. For instance, the {\it primitive cubic lattice} (the lattice $\mathbb{Z}^3$ with the standard network structure) is symmetric, but not strongly isotropic.}

\item It is a web of congruent {\it decagonal rings} (minimal circuits of length 10 in the graph-theoretical sense; Fig.~\!\ref{deca}). There are 15 decagonal rings passing through each vertex (Fig.~\!\ref{15}; see Remark \ref{rem:ccc}(2)).  

\begin{figure}[htbp]
\vspace{-0.4cm}
\hspace{3.9cm}
\includegraphics[width=0.4\linewidth]{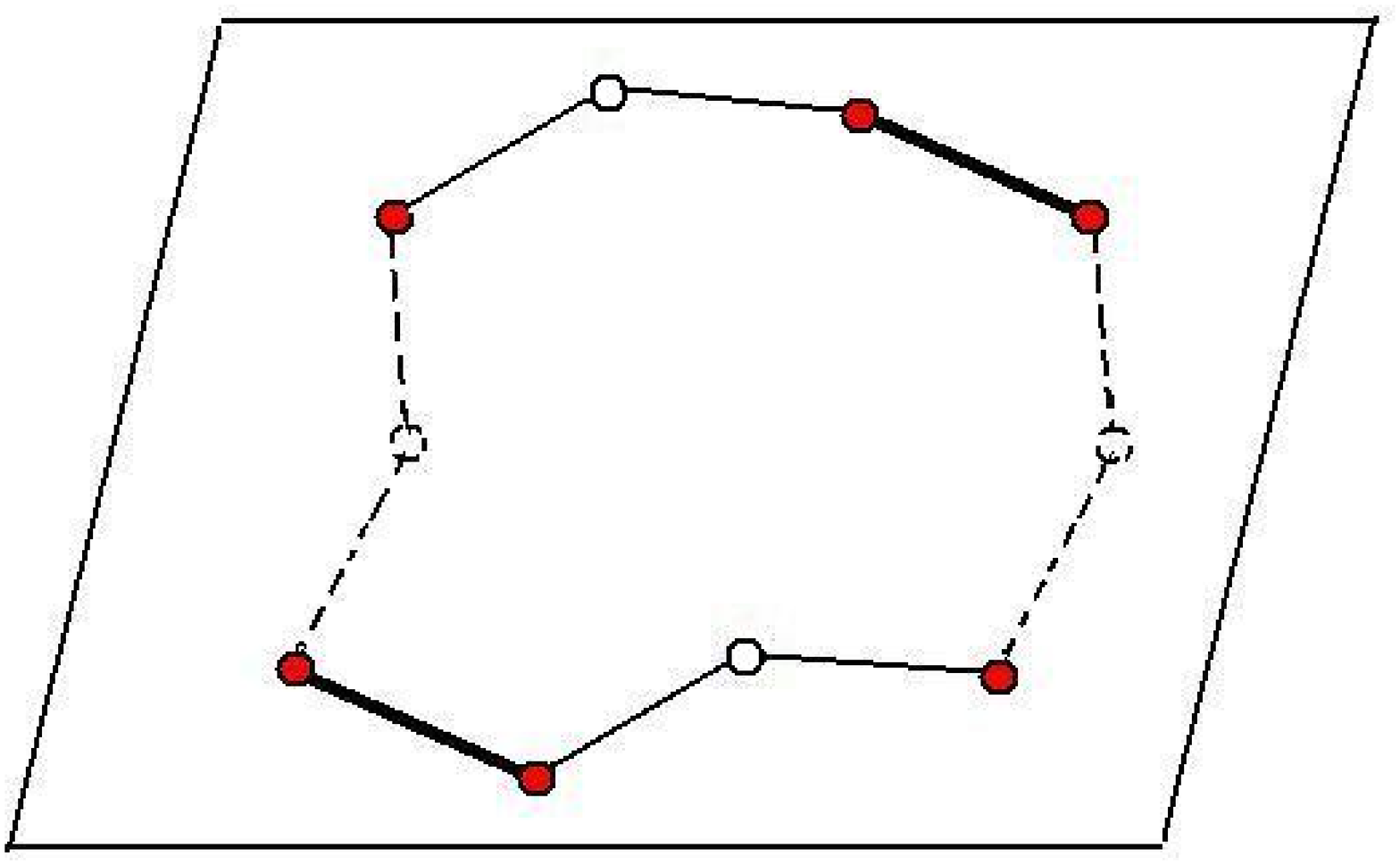}
\vspace{-0.6cm}
\caption{Decagonal ring}\label{deca}
\end{figure}

\item It is characterized by the {\it minimizing property} for a certain energy functional.\footnote{Regarding a crystal as a system of harmonic oscillators, we  may define ``energy per unit cell". We should note that it is different from energy in the physical sense (cf. \cite{shu}). See \cite{kotani2}, \cite{kotani1}, \cite{su6} for the detail. Laves' graph also minimizes the total edge length normalized by volume (\cite{alex}).}

\item There exists a Cartesian coordinate system such that each vertex has an integral coordinate.

\item It has {\it chirality}; that is, it is non-superposable on its mirror image.

\end{enumerate}

As the first two properties tell us, Laves' graph has very big symmetry. Among them, Property 1 implies that the automorphism group is identified with a 3-dimensional crystallographic group, which, as usual in crystallography, describes macroscopic symmetry of the crystal. Its point group turns out to be the {\it octahedral group}, the group of rotations fixing the regular octahedron. Property 3 is most peculiar among all. It is truly amazing that such a big number of congruent decagonal rings may gather at each vertex in the same manner.

\begin{figure}[htbp]
\vspace{-0.5cm}
\hspace{1.2cm}
\includegraphics[width=0.8\linewidth]{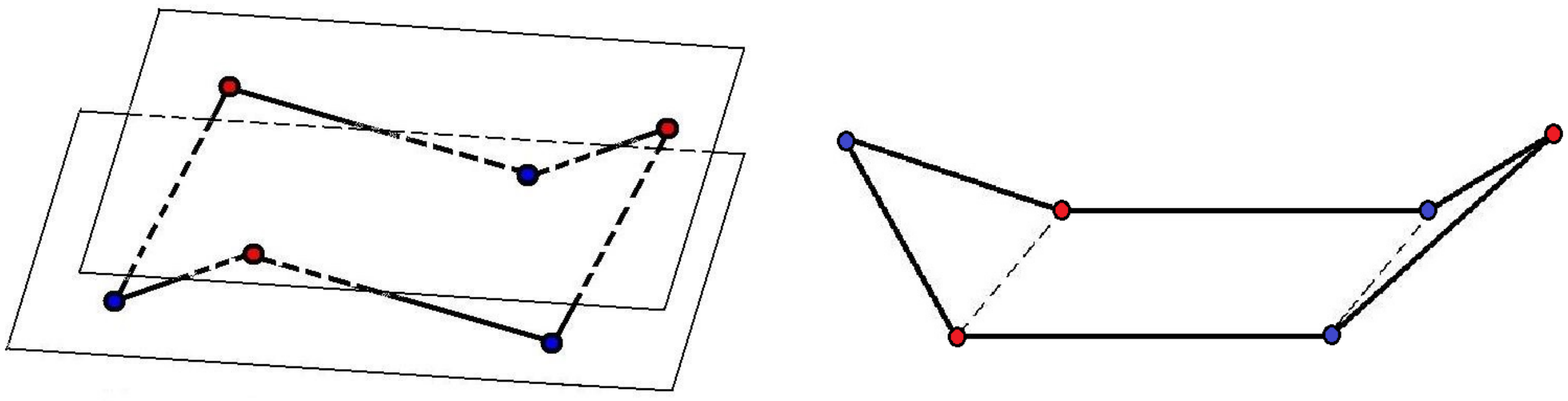}
\vspace{-4.2cm}
\caption{chair and boat conformations}\label{cb}
\end{figure}

To justify the name ``diamond twin" for Laves' graph, we look at the structure of the (cubic) diamond, a real crystal with very big microscopic symmetry.\footnote{Besides carbon, there are several elements (say, silicon (Si), germanium (Ge), tin (Sn)) which adopt the same structure as the diamond. By abuse of language,  the term ``diamond" is used to express its network structure throughout, not to stand for the diamond as an actual crystal.} The diamond structure is the standard realization of the maximal abelian covering graph of the {\it dipole graph} with two vertices joined by 4 parallel edges, and is a web of congruent hexagonal rings (called {\it chair conformation}; see the left of Fig.~\!\ref{cb}). Moreover, it has Properties 1, 2, 4, 5. The point group is the {\it full octahedral group}, the group consisting of all orthogonal transformations fixing the regular octahedron.\footnote{The full octahedral group, symbolically denoted by ${\rm O}_{\rm h}$ in {\it Sch\"{o}nflies notation}, is isomorphic to the direct product $\mathcal{S}_4\times \mathbb{Z}_2$, where $\mathcal{S}_4$ is the symmetry group of four objects. On the other hand, the octahedral group---isomorphic to $\mathcal{S}_4$ and symbolically denoted by ${\rm O}$---is the group consisting of all rotations acting on the octahedron.} 
The number of hexagonal rings passing through each vertex is 12, which is less than 15 but still a big number.\footnote{Another example of a crystal structure having a big number of congruent rings is the primitive cubic lattice, which is a web of four-membered rings forming squares. The number of such rings passing through each vertex is 12.}

What is more, a crystal having Properties 1, 2, 4 must be either diamond or Laves' graph (or its mirror image because of chirality) (\cite{su1}). In this sense, Laves' graph and the diamond structure are very kinfolk as mathematical objects (a difference is that the diamond structure has no chirality).

\begin{figure}[htbp]
\vspace{-0.1cm}
\begin{center}
\includegraphics[width=0.8\linewidth]{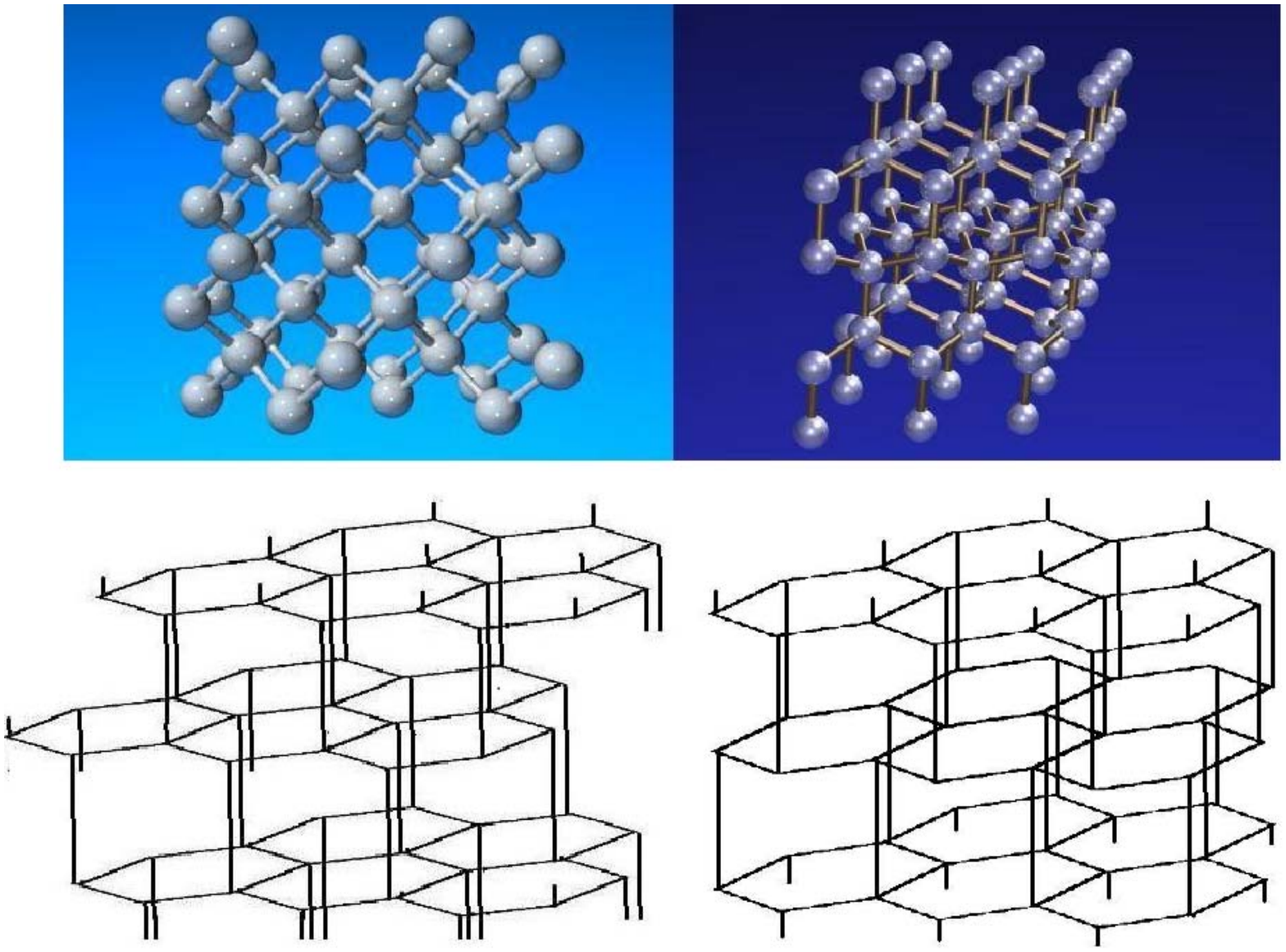}
\end{center}
\vspace{-0.8cm}
\caption{Diamond and Lonsdaleit. From https://www.webelements.com}\label{DL}
\end{figure}

Here is one remark in order. The diamond twin does not belong to the family of crystal structures of the so-called {\it diamond polytypes} (or {\it diamond cousins} in plain language) such as {\it Lonsdaleite} (named in honor of Kathleen Lonsdale and also called {\it Hexagonal diamond}, a rare stone of pure carbon discovered at Meteor Crater, Arizona, in 1967). Incidentally, the network structure of Lonsdaleite---having much less symmetry than the diamond---is a web of two types of congruent hexagonal rings; one being in the chair conformation, and another being in the {\it boat conformation} (the right of Fig.~\!\ref{cb}). A shape-similarity between the structures of diamond and Lonsdaleite is brought out when looking at the graphite-like realizations of those structures (see the lower figures in Fig.~\!\ref{DL}),\footnote{In general, a graphite-like realization means a crystal net with a layered structure that consists of 2-dimensional crystal nets placed in  horizontal sheets.} therefore if we stick to the apparent shape, not to symmetry enshrined inward, it might be appropriate to call Lonsdaleite the twin of diamond, but we do not adopt this view.

This note is a spin-off of the previous articles \cite{su1} and \cite{sunada1}. The aim is to reinforce the reason for the appellation ``diamond twin," taking a look at the {\it period lattices} and the {\it building blocks} of the diamond crystal and its twin. With these instruments, we first make a down-to-earth observation about mutual relations between the two crystal structures. For instance, it is observed that the building blocks of the diamond and its twin are interwoven with each other via their period lattices. What is more remarkable is, as seen in the last section, that each of the building blocks generates an ``orthogonally symmetric lattice," an interesting notion in its own right and a natural generalization of {\it irreducible root lattices} (the concept stemming from the theory of Lie groups and Lie algebras). This notion, combined with an exhaustive search of finite subgroups of the 3-dimensional orthogonal group ${\rm O}(3)$, makes the diamond and its twin very distinct among all crystal structures. The argument taken up 
is elementary in nature, but requires a careful enumeration of many possible cases.

In a nutshell, what we attempt to do is to provide a specific example that concretizes Kepler's famous statement ``At ubi materia, ibi Geometria" (where there is matter, there is geometry).\footnote{J. Kepler,
{\it De Fundamentis Astrologiae Certioribus} (Concerning the More Certain Fundamentals of Astrology), 1601. In this regard, it is worth recalling that geometry in ancient Greece---especially the classification of regular polyhedra that is regarded as the culmination of Euclid's {\it Elements}---had its source in the curiosity to the shapes of crystals, as it is often said.} 

\section{Building blocks of the diamond and its twin} 
We shall very briefly review some materials in graph theory and elementary algebraic topology to explain the notion of building blocks associated with crystal nets. Although the objects of our concern is 3-dimensional crystals,  we deal here with crystals of general dimension. See \cite{su3}, \cite{su6} for the details.

A {\it graph}\index{graph} is represented by an ordered pair $X = (V,E)$ of the set of {\it vertices} $V$ and the set of all {\it directed edges} $E$. For a  directed edge $e$, we denote by ${\it o}(e)$ the {\it origin}, and by ${\it t}(e)$ the {\it terminus}. The inversed edge of $e$ is denoted by $\overline{e}$. The set of directed edges with origin $x \in V$ is denoted by $E_x$; i.e., $
E_x = \{ e \in E |~\!o(e)=x \}$.

The net associated with a $d$-dimensional crystal is not just an infinite graph realized in $\mathbb{R}^d$, but a graph with a translational action of a {\it lattice} (called ``period  lattice")\footnote{A lattice is a discrete subgroup of the additive group $\mathbb{R}^d$ of maximal rank. More specifically, $L\subset \mathbb{R}^d$ is a lattice if there exists a basis $\{{\bf a}_1,\ldots,{\bf a}_d\}$ of $\mathbb{R}^d$ such that $L=\{k_1{\bf a}_1+\cdots+k_d{\bf a}_d|~\!k_i\in \mathbb{Z}~(i=1,\ldots,d)\}$. $\{{\bf a}_1,\ldots,{\bf a}_d\}$ is called a $\mathbb{Z}$-basis of $L$.} which becomes a finite graph when factored out.\footnote{This fact was pointed out by crystallographers (\cite{chung}). Mathematically, a crystal net as an abstract graph is an infinite-fold abelian covering graph of a finite graph.} The finite graph obtained by factoring out is called the {\it quotient graph}. 

We let $X=(V,E)$ be the abstract graph associated with a $d$-dimensional crystal net with a period lattice $L$, and let $X_0=(V_0,E_0)$ be the quotient graph. 
We assign a vector ${\bf v}(e)$ to each directed edge $e$ in $X_0$ as follows. Choose a direct edge $e'$ in $X$ which corresponds to $e$. In the crystal net, $e'$ is a directed line segment, so $e'$ yields a vector ${\bf v}(e)\in \mathbb{R}^d$, which, as easily checked, does not depend upon the choice of $e'$. Obviously ${\bf v}(\overline{e})=-{\bf v}(e)$. Hence, a building block is nothing but a $1$-{\it cochain} on $X_0$ with values in $\mathbb{R}^d$.

\begin{figure}[htbp]
\vspace{-0.7cm}
\begin{center}
\includegraphics[width=0.8\linewidth]{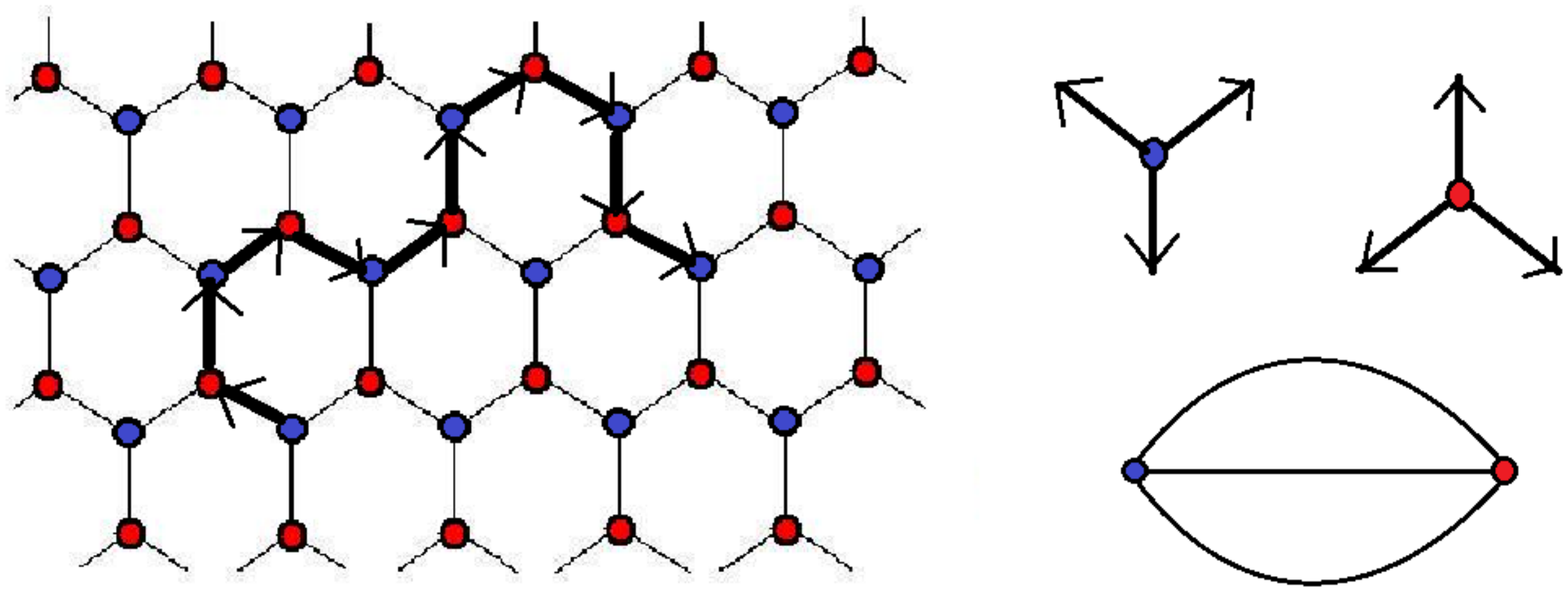}
\end{center}
\vspace{-3.7cm}
\caption{Honeycomb lattice}\label{honey}
\end{figure}

Put ${\bf E}_{x}={\bf v}(E_{0x})$. The system of vectors
$
\{{\bf E}_x\}_{x\in V_0}
$ completely determines the original crystal net and its period lattice. In fact, we obtain the original crystal net by summing up vectors ${\bf v}(e_i)$ for all paths $(e_1,\ldots,e_n)$ ($e_i\in E_0$) on $X_0$ which start from a reference vertex.\footnote{A path means a sequence of edges $(e_1,\ldots,e_n)$ with $t(e_i)=o(e_{i+1})$ $(i=1,\ldots,n-1)$. 
An example is depicted in Fig.~\!\ref{honey} by a sequence of arrows leading one after another.} The period lattice $L$ turns out to be the image $\widehat{{\bf v}}(H_1(X_0,\mathbb{Z}))$ of the homomorphism 
$
\widehat{{\bf v}}:H_1(X_0,\mathbb{Z})\longrightarrow \mathbb{R}^d\quad (d=2,3)
$ 
defined by 
$\widehat{{\bf v}}\big(\sum_{e\in E_0}a_ee\big)=\sum_{e\in E_0}a_e{\bf v}(e)$, where $H_1(X_0,\mathbb{Z})$ is the 1st homology group of $X_0$. Thus, the system $\{{\bf E}_x\}_{x\in V_0}$  deserves to be called the {\it building block} of the crystal.

For later purposes, we make a small remark: In the case that $X$ is the maximal abelian covering graph of $X_0$ (thus $d={\rm rank}~\!H_1(X_0,\mathbb{Z})$), one may take closed paths $c_1,\ldots,c_d$ in $X_0$ such that $\widehat{{\bf v}}([c_1]),\ldots,\widehat{{\bf v}}([c_d])$ comprise a basis of the period lattice, where $[c_i]$ is the homology class represented by $c_i$.

Figure \ref{honey} illustrates how the honeycomb lattice is obtained from the building block and the quotient graph.

We shall now describe explicitly the building blocks for the diamond and its twin in the coordinate space $\mathbb{R}^3$. To facilitate understanding of the configurations of vectors in the building blocks, we employ the cube $Q$ as an auxiliary figure whose vertices are $(1,1,1), (1,-1,1),$ $(-1,1,1),$ $(-1,-1,1),$ $(1,1,-1),$ $(1,-1,-1),$ $(-1,1,-1),$ $(-1,-1,-1)$.
 
\bigskip

\noindent{\bf (I)~ The building block for the diamond twin}

\begin{figure}[htbp]
\vspace{-1cm}
\hspace{-0.8cm}
\includegraphics[width=1.1\linewidth]{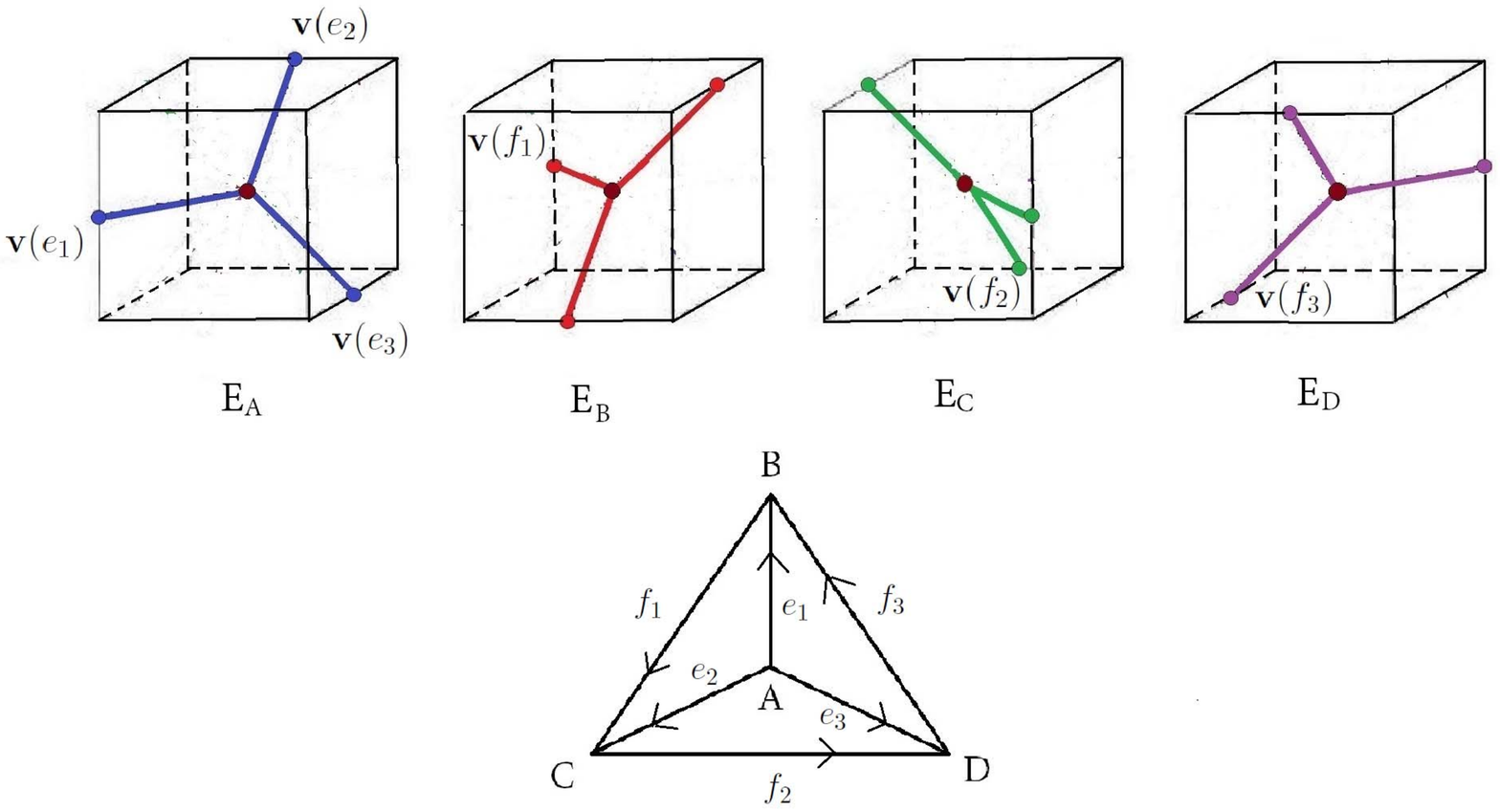}
\vspace{-2.8cm}
\caption{The building block and the quotient graph for the diamond twin}\label{building1}
\end{figure}

Let $X_0$ be the complete graph $K_4$ with vertices $A,B,C,D$ (see the lower figure in Fig.~\!\ref{building1}), and put
\begin{eqnarray*}
&&{\bf E}_A=\{{}^t(0,1,1), {}^t(-1,-1,0), {}^t(1,0,-1)\},\\
&&{\bf E}_B=\{{}^t(1,0,1), {}^t(-1,1,0), {}^t(0,-1,-1)\},\\
&&{\bf E}_C=\{{}^t(-1,0,1), {}^t(0,1,-1), {}^t(1,-1,0)\},\\
&&{\bf E}_D=\{{}^t(0,-1,1), {}^t(-1,0,-1), {}^t(1,1,0)\},
\end{eqnarray*}
where vectors in $\mathbb{R}^3$ are represented by column vectors. The system 
$\{{\bf E}_A,{\bf E}_B,$ ${\bf E}_C, {\bf E}_D\}$ forms the building block for the diamond twin (see the upper diagrams of Fig.~\!\ref{building1}). 
Note that each of ${\bf E}_A,$ ${\bf E}_B,$ ${\bf E}_C, {\bf E}_D$ comprises an equilateral triangle in a plane with barycenter ${\bf o}=(0,0,0)$. More specifically, the planes containing ${\bf E}_A, {\bf E}_B, {\bf E}_C,{\bf E}_D$ are orthogonal to the vectors ${\bf a}={}^t(1,-1,1), {\bf b}={}^t(1,1,-1), {\bf c}={}^t(-1,-1,$ $-1), {\bf d}={}^t(-1,1,1)$, respectively, from which we see that the dihedral angle of any two planes is $\arccos 1/3\fallingdotseq 70.53^{\circ}$.

In view of Fig.~\!\ref{building1}, we have ${\bf v}(e_1)={}^t(-1,-1,0), {\bf v}(e_2)={}^t(0,1,1), {\bf v}(e_3)={}^t(1,0,-1)$ and ${\bf v}(f_1)={}^t(-1,1,0), {\bf v}(f_2)={}^t(0,1,-1), {\bf v}(f_3)={}^t(-1,0,-1)$.
As a $\mathbb{Z}$-basis of $H_1(X_0,\mathbb{Z})$, one can take $[c_1],[c_2],[c_3]$ where $c_1=(e_2,f_1,\overline{e}_3)$,~ $c_2=(e_3,f_2, \overline{e}_1)$,~ $c_3=(e_1,f_3,\overline{e}_2)$. We then have
\begin{eqnarray*}
&&\widehat{{\bf v}}([c_1])={}^t(0,1,1)+{}^t(-1,1,0)+{}^t(-1,0,1)=2\cdot{}^t(-1,1,1),\\
&&\widehat{{\bf v}}([c_2])={}^t(1,0,-1)+{}^t(0,1,-1)+{}^t(1,1,0)=2\cdot{}^t(1,1,-1),\\
&&\widehat{{\bf v}}([c_3])={}^t(-1,-1,0)+{}^t(-1,0,-1)+{}^t(0,-1,-1)=2\cdot{}^t(-1,-1,-1),\end{eqnarray*}
which comprise a $\mathbb{Z}$-basis of the period lattice of the diamond twin. 

We let $L_{\mathcal{DT}}$ be the lattice with $\mathbb{Z}$-basis ${}^t(-1,1,1)$, ${}^t(1,1,-1)$, ${}^t(-1,-1,-1)$ (hence $2L_{\mathcal{DT}}$ is the period lattice of the diamond twin). It is checked that 
$$
L_{\mathcal{DT}}=\{(x_1,x_2,x_3)\in \mathbb{Z}^3|~\!x_1+x_2,~x_2+x_3,~x_3+x_1~\text{are even}\}.
$$
Indeed, if we write $(x_1,x_2,x_3)=k_1(-1,1,1)+k_2(1,1,-1)+k_3(-1,-1,-1)$, then
$x_1+x_2=2(k_2-k_3),$ $x_2+x_3=-2k_3$, $x_3+x_1=2(-k_1-k_3)$, and
$k_1=\frac{1}{2}(x_1+x_2)-x_1$, $k_2=\frac{1}{2}(x_1+x_3)-x_3$, $k_3=-\frac{1}{2}(x_2+x_3)$.

\begin{figure}[htbp]
\vspace{-0.4cm}
\begin{center}
\includegraphics[width=.5\linewidth]{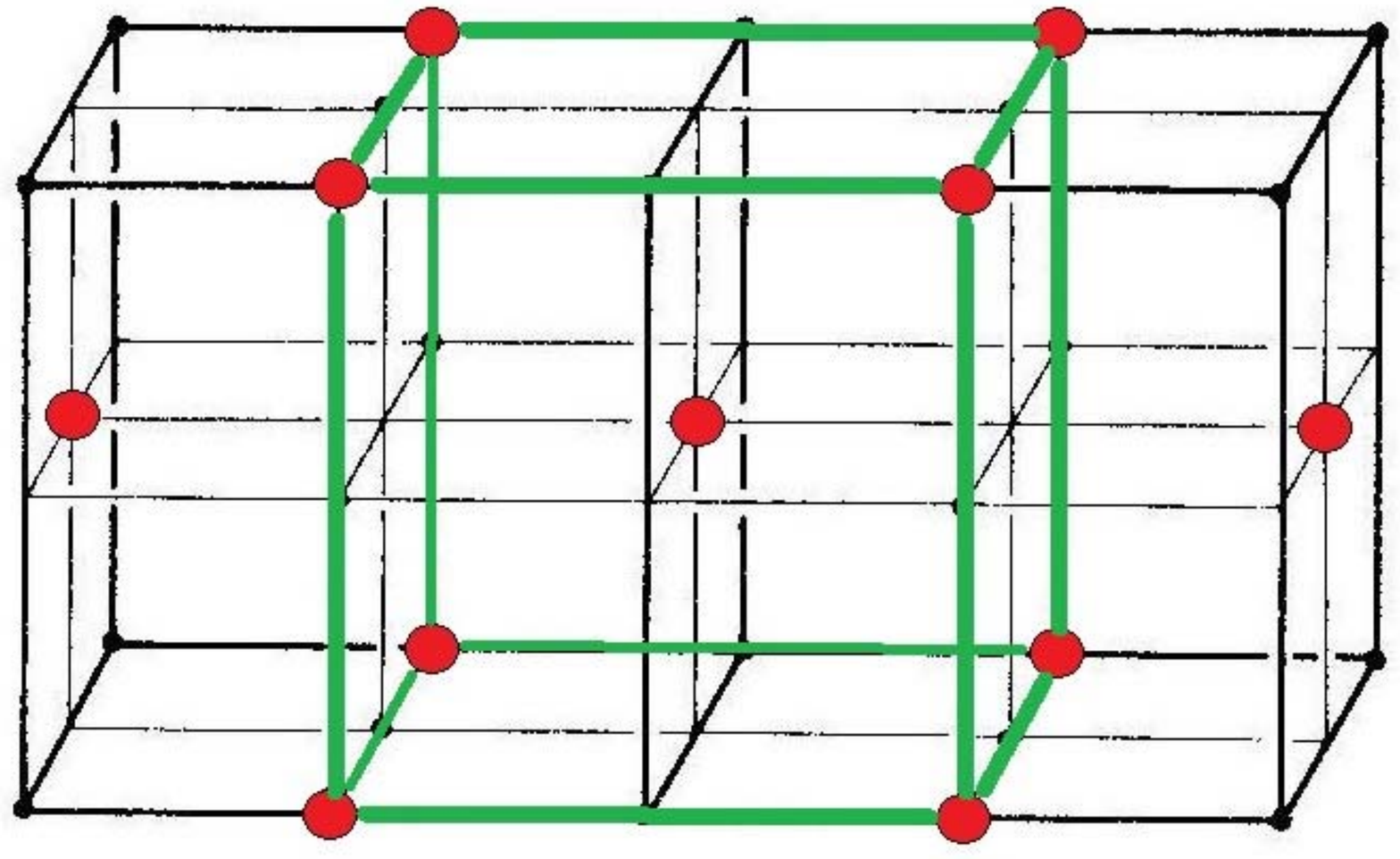}
\end{center}
\vspace{-1.2cm}
\caption{Body-centered cubic lattice}\label{body}
\end{figure}

The lattice $L_{\mathcal{DT}}$ is what is called the {\it body-centered cubic lattice} in crystallography (look at the cube in Fig.~\!\ref{body} depicted by the bold lines).

\medskip
\noindent{\bf (II)~ The building block for the diamond}

\begin{figure}[htbp]
\vspace{-0.8cm}
\hspace{1.9cm}
\includegraphics[width=.68\linewidth]{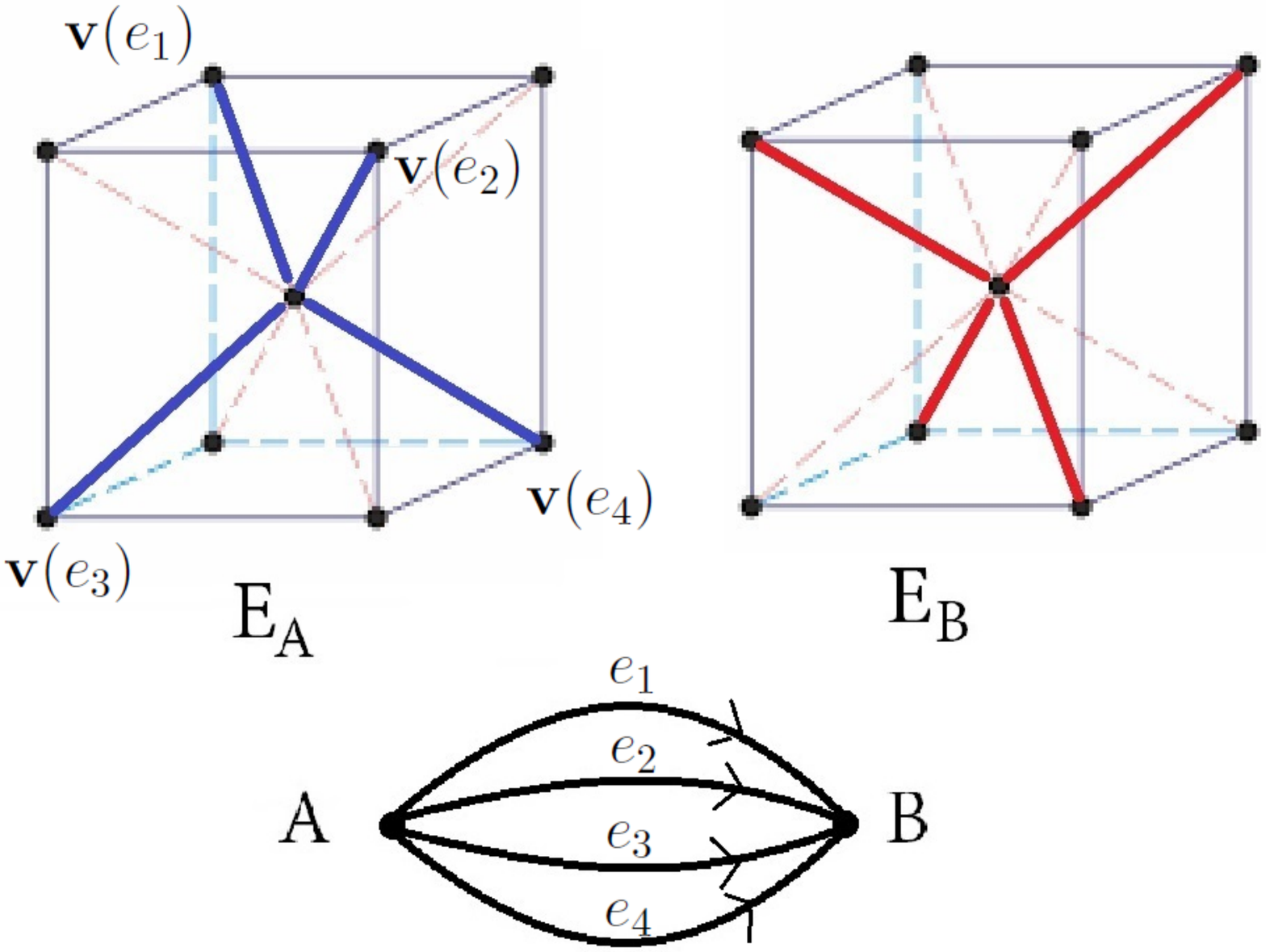}
\vspace{-0.6cm}
\caption{The building block and the quotient graph for the diamond}\label{building2}
\end{figure}

Let $X_0$ be the graph with two vertices $A, B$ joined by 4 parallel edges $e_1,e_2,e_3,e_4$ (see the lower figure in Fig.~\!\ref{building2}),
and put
\begin{eqnarray*}
&&{\bf E}_A=\{{}^t(-1,1,1), {}^t(1,-1,1), {}^t(-1,-1,-1), {}^t(1,1,-1)\},\\
&&{\bf E}_B=\{{}^t(1,1,1), {}^t(-1,-1,1), {}^t(-1,1,-1), {}^t(1,-1,-1)\}=-{\bf E}_A.
\end{eqnarray*}
The system 
$\{{\bf E}_A,{\bf E}_B\}$ forms the building block for the diamond (see the upper diagrams of Fig.~\!\ref{building2}). Note that each of ${\bf E}_A$ and ${\bf E}_B$ comprises a regular tetrahedron with barycenter ${\bf o}$.\footnote{This description tells us that the honeycomb lattice (the {\it graphene} as a real matter) is regarded as a 2-dimensional analogue of the diamond crystal.}

As shown in Fig.~\!\ref{building2}, we put ${\bf v}(e_1)={}^t(-1,1,1)$, ${\bf v}(e_2)={}^t(1,-1,1)$, ${\bf v}(e_3)={}^t(-1,-1,-1)$, ${\bf v}(e_4)={}^t(1,1,-1)$.
As a $\mathbb{Z}$-basis of $H_1(X_0,\mathbb{Z})$, one can take $[c_1],$ $[c_2],$ $[c_3]$ where $c_1=(e_1,\overline{e}_2)$,~ $c_2=(e_2, \overline{e}_3)$,~ $c_3=(e_3,\overline{e}_4)$. We then have
\begin{eqnarray*}
&&\widehat{{\bf v}}([c_1])={}^t(-1,1,1)+{}^t(-1,1,-1)=2\cdot{}^t(-1,1,0),\\
&&\widehat{{\bf v}}([c_2])={}^t(1,-1,1)+{}^t(1,1,1))=2\cdot{}^t(1,0,1),\\
&&\widehat{{\bf v}}([c_3])={}^t(-1,-1,-1)+{}^t(-1,-1,1)=2\cdot{}^t(-1,-1,0),\end{eqnarray*}
which comprise a $\mathbb{Z}$-basis of the period lattice of the diamond.

\begin{figure}[htbp]
\vspace{-0.4cm}
\begin{center}
\includegraphics[width=.48\linewidth]{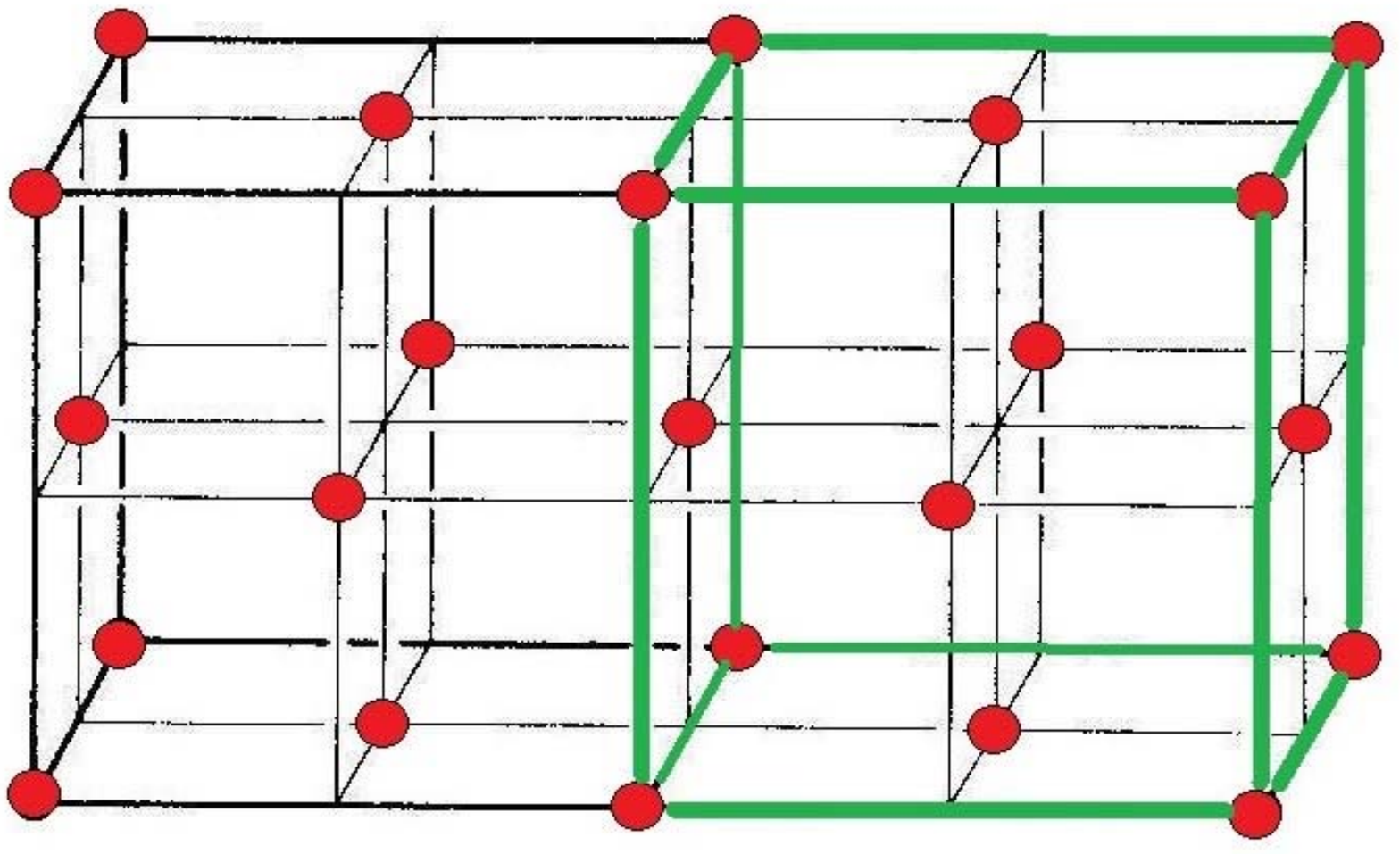}
\end{center}
\vspace{-1.2cm}
\caption{Face-centered cubic lattice}\label{face}
\end{figure}

We let $L_{\mathcal{D}}$ be the lattice with $\mathbb{Z}$-basis ${}^t(-1,1,0)$, ${}^t(1,0,1)$, ${}^t(-1,-1,0)$ (hence $2L_{\mathcal{D}}$ is the period lattice of diamond). By the same method as described in (I), it is checked that 
$$
L_{\mathcal{D}}=
\{(x_1,x_2,x_3)\in \mathbb{Z}^3|~\!x_1+x_2+x_3~\text{is even}\}.
$$
This is what is called the {\it face-centered cubic lattice} in crystallography (look at the cube in Fig.~\!\ref{face} depicted by the bold lines).

\medskip

\begin{figure}[htbp]
\vspace{-0.4cm}
\begin{center}
\includegraphics[width=.61\linewidth]{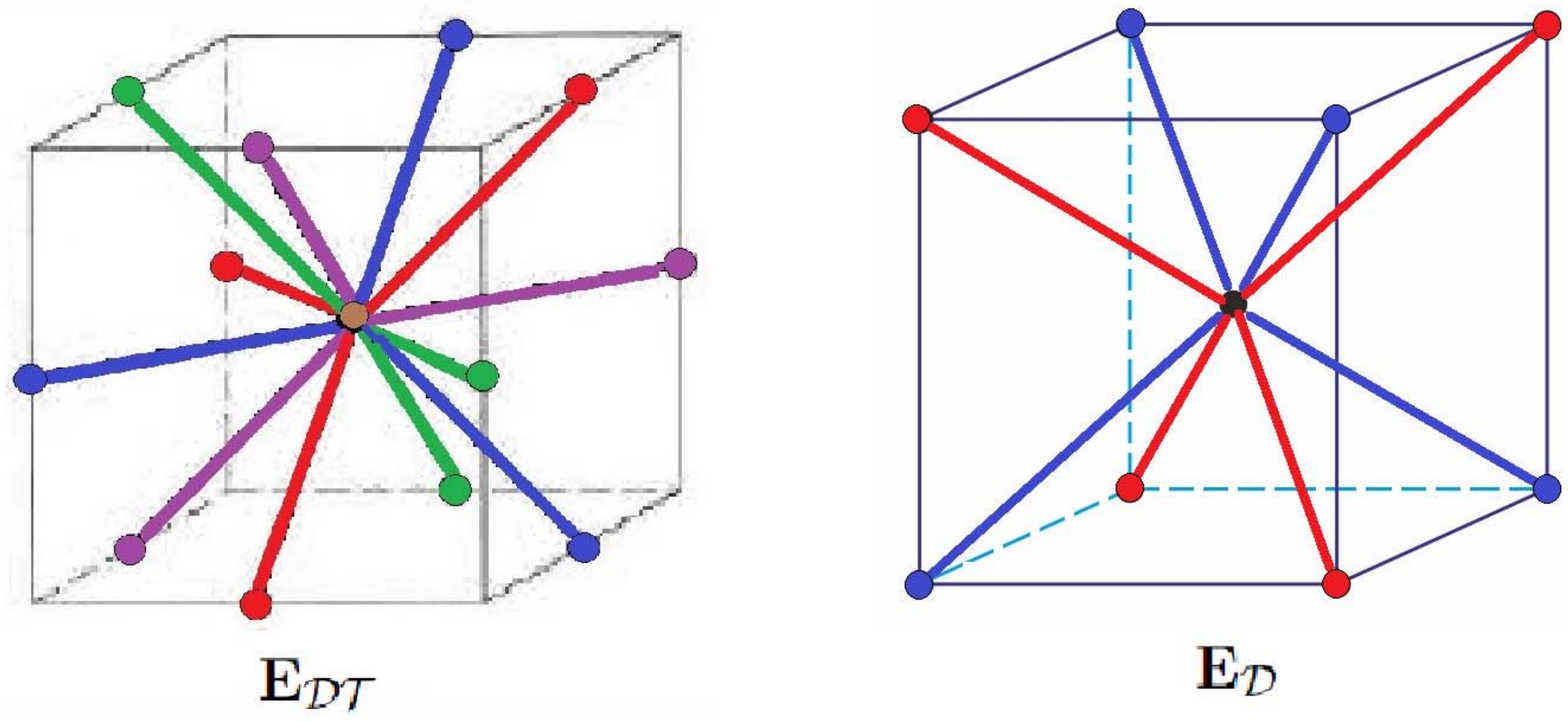}
\end{center}
\vspace{-2.7cm}
\caption{${\bf E}_{\mathcal{DT}}$ and ${\bf E}_{\mathcal{D}}$}\label{union}
\end{figure}

So far, no relationship can be found between the diamond and its twin.
A mutual relation between them, though not a big deal, is observed only after considering the union ${\bf E}_{\mathcal{DT}}={\bf E}_A\cup {\bf E}_B\cup {\bf E}_C\cup {\bf E}_D$ for the diamond twin and the union ${\bf E}_{\mathcal{D}}={\bf E}_A\cup {\bf E}_B$ for the diamond (Fig.~\!\ref{union}). Perhaps both systems of vectors may be familiar to the reader. For instance, ${\bf E}_{\mathcal{DT}}$ is nothing but the irreducible root system $A_3$ (see the next section).

The system ${\bf E}_{\mathcal{DT}}$ for the diamond twin generates the lattice $L_{\mathcal{D}}$ since ${\bf E}_{\mathcal{DT}}\subset L_{\mathcal{D}}$ and contains a basis of $L_{\mathcal{D}}$. (Remember 
 that $2L_{\mathcal{D}}$ is the period lattice for the diamond, not for the diamond twin!). 
On the other hand, the system ${\bf E}_{\mathcal{D}}$ for the diamond generates the lattice $L_{\mathcal{DT}}$ since ${\bf E}_{\mathcal{D}}\subset L_{\mathcal{DT}}$ and contains a basis of $L_{\mathcal{DT}}$. Hence, passing from the building blocks to the period lattices, the role exchange between the diamond and its twin takes place. 

Furthermore, if we denote by $L^*$ the dual (reciprocal) lattice of a lattice $L\subset \mathbb{R}^d$ in general,\footnote{$L^*=\{{\bf x}\in \mathbb{R}^d|~\!\langle {\bf x},{\bf y}\rangle\in \mathbb{Z}~\text{for every ${\bf y}\in L$}\}$.} we have
$$
L_{\mathcal{DT}}{}^*=\frac{1}{2}L_{\mathcal{D}}, ~~L_{\mathcal{D}}{}^*=\frac{1}{2}L_{\mathcal{DT}}
$$ 
since the dual of the $\mathbb{Z}$-basis $\{{}^t(-1, 1, 1), {}^t(1,1,-1), {}^t(-1,-1,-1)\}$ of $L_{\mathcal{DT}}$ is $\{{}^t(-1/2,1/2,0),$ ${}^t(0,1/2,-1/2), {}^t(-1/2,0,-1/2)\}$. Here we should note that $\{{}^t(-1,1,0),$ ${}^t(0,1,-1), {}^t(-1,0,-1)\}$ is a $\mathbb{Z}$-basis of $L_{\mathcal{D}}$.

What is more peculiar than the facts obtained by the simple observations above is that $L_{\mathcal{D}}$ and $L_{\mathcal{DT}}$ are {\it orthogonally symmetric lattices}, the notion introduced in the next section which distinguishes the diamond and its twin from all other crystal structures (with one exception).

\begin{remark}\label{rem:ccc} {\rm 
(1)~ Just for the sake of completeness, let us give a precise description of the network structure of the diamond twin.
Put $p_1=(0,1,1),~p_2=(1,0,-1),~p_3=(-1,-1,0)$.
The set of vertices ${\bf V}$ of the diamond twin is decomposed as 
$$
{\bf V}={\bf V}_0\cup {\bf V}_1\cup {\bf V}_2\cup {\bf V}_3, 
$$
where
\begin{eqnarray*}
&&{\bf V}_0=2L_{\mathcal{DT}},~~
{\bf V}_1=p_1+2L_{\mathcal{DT}},\\
&&{\bf V}_2=p_2+2L_{\mathcal{DT}},~~
{\bf V}_3=p_3+2L_{\mathcal{DT}}.
\end{eqnarray*}
Note that the sets ${\bf V}_0, {\bf V}_1, {\bf V}_2, {\bf V}_3$ correspond to the vertices $A,B,C,D$ in the quotient graph $K_4$, respectively. 

The incidence relation of vertices is described as follows:
\begin{enumerate}
\item $\alpha\in {\bf V}_0$ and $p_i+\beta\in {\bf V}_i$ ($i=1,2,3$) are joined by an edge if and only if $\alpha=\beta$.
\item $p_1+\alpha\in {\bf V}_1$ and $p_2+\beta\in {\bf V}_2$ are joined by an edge if and only if $\beta=\alpha+(-2,2,2)$.
\item $p_2+\alpha\in {\bf V}_2$ and $p_3+\beta\in {\bf V}_3$ are joined by an edge if and only if $\beta=\alpha+(2,2,-2)$.
\item $p_3+\alpha\in {\bf V}_3$ and $p_1+\beta\in {\bf V}_1$ are joined by an edge if and only if $\beta=\alpha+(-2,-2,-2)$.
\end{enumerate}

\medskip

(2)~ 
As mentioned in Introduction, the most peculiar feature of the diamond twin is that the number of decagonal rings passing through each vertex is $15$. However, it is not easy to infer this fact 
from the CG image or an available network model because of its intricate space-pattern. Here let us give an explicit description of those decagonal rings in terms of paths in the quotient graph.
. 

Consider the graph $K_4$ with labeled edges as Fig.~\!\ref{building1}. The decagonal rings passing through a vertex 
are obtained from the following closed paths in $K_4$ all of which are homologous to zero.

\medskip

(I) \hspace{0.05cm} \qquad  $(e_1,f_3,\overline{e_2},e_3,f_2,\overline{e_1},e_2,\overline{f_3},\overline{f_2},\overline{e_3})$

(II) \qquad  $(e_1,\overline{f_2},\overline{e_3},e_2,\overline{f_3},\overline{e_1},e_3,f_2,f_3,\overline{e_2})$

(III)\qquad  $(e_2,\overline{f_3},\overline{e_1},e_3,f_2,f_3,\overline{e_2},e_1,\overline{f_2},\overline{e_3})$

(IV)\qquad  $(e_1,f_3,\overline{e_2},e_3,\overline{f_1},\overline{f_3},\overline{e_1},e_2,f_1,\overline{e_3})$

(V)\hspace{0.1cm}\qquad  $(e_1,f_3,f_1,\overline{e_3},e_2,\overline{f_3},\overline{e_1},e_3,\overline{f_1},\overline{e_2})$

(VI)\qquad  $(e_2,\overline{f_3},\overline{e_1},e_3,\overline{f_1},\overline{e_2},e_1,f_3,f_1,\overline{e_3})$

(VII)\hspace{0.18cm}\quad  $(e_1,\overline{f_2},\overline{e_3},e_2,f_1,f_2,\overline{e_1},e_3,\overline{f_1},\overline{e_2})$

(VIII)\hspace{0.06cm}\quad  $(e_1,\overline{f_2},\overline{f_1},\overline{e_2},e_3,f_2,\overline{e_1},e_2,f_1,\overline{e_3})$

(IX)\qquad  $(e_2,f_1,f_2,\overline{e_1},e_3,\overline{f_1},\overline{e_2},e_1,\overline{f_2},\overline{e_3})$

(X)\hspace{0.1cm}\qquad  $(e_1,f_3,f_1,f_2,\overline{e_1},e_3,\overline{f_1},\overline{f_3},\overline{f_2},\overline{e_3})$
 
(XI)\qquad  $(e_1,\overline{f_2},\overline{f_1},\overline{f_3},\overline{e_1},e_3,f_2,f_3,f_1,\overline{e_3})$

(XII)\hspace{0.18cm}\quad  $(e_1,f_3,f_1,f_2,\overline{e_1},e_2,\overline{f_3},\overline{f_2},\overline{f_1},\overline{e_2})$

(XIII)\hspace{0.05cm}\quad $(e_1,\overline{f_2},\overline{f_1},\overline{f_3},\overline{e_1},e_2,f_1,f_2,f_3,\overline{e_2})$

(XIV)\hspace{0.05cm}\quad  $(e_2,f_1,f_2,f_3,\overline{e_2},e_3,\overline{f_1},\overline{f_3},\overline{f_2},\overline{e_3})$

(XV)\hspace{0.18cm}\quad  $(e_2,\overline{f_3},\overline{f_2},\overline{f_1},\overline{e_2},e_3,f_2,f_3,f_1,\overline{e_3})$

}

\end{remark}

\section{Orthogonally symmetric lattices}
\begin{figure}[htbp]
\vspace{-0.6cm}
\begin{center}
\includegraphics[width=.85\linewidth]{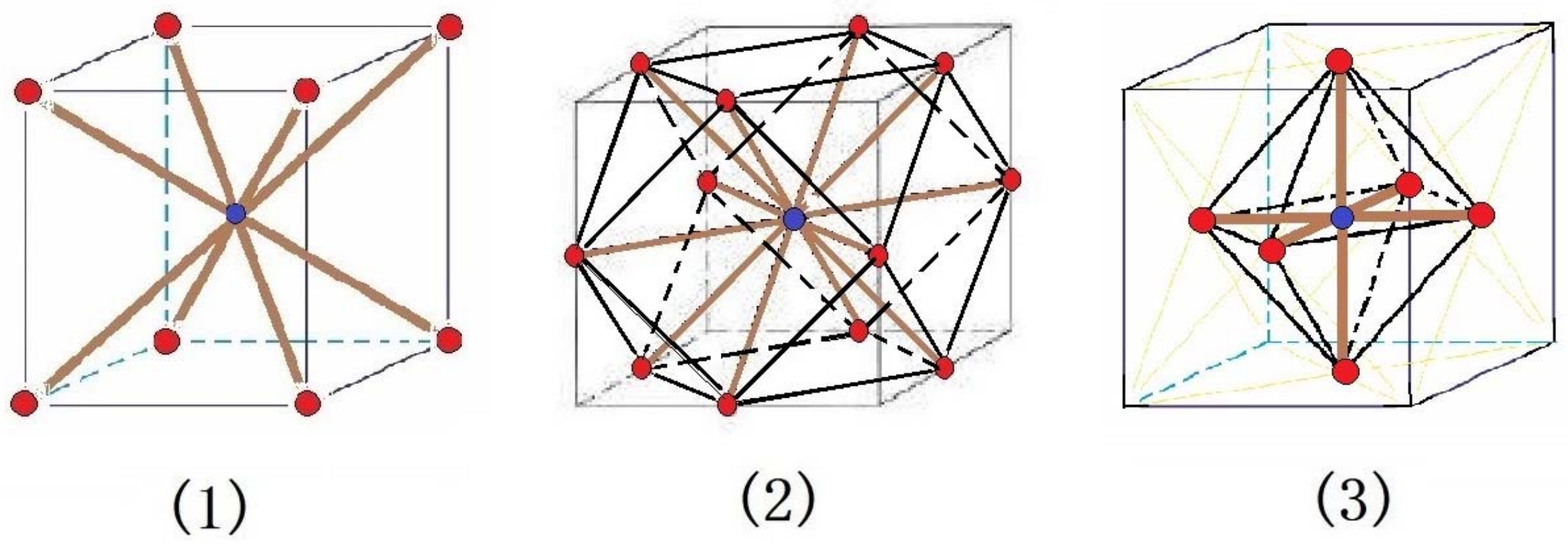}
\end{center}
\vspace{-4.2cm}
\caption{}\label{octa3}
\end{figure}

We denote by $S_{r}({\bf a})$ the sphere in $\mathbb{R}^d$ of radius $r$, centered at ${\bf a}\in \mathbb{R}^d$. 
Given a general lattice $L$ in $\mathbb{R}^d$, we put $\alpha(L):=\min_{{\bf y}\neq{\bf 0}\in L}\|{\bf y}\|$, and 
\begin{eqnarray*}
&&K(L):=\{{\bf x}\in L|~\!\|{\bf x}\|=\alpha(L)\},\\
&&G(L):=\{g\in {\rm O}(d)|~\!g(L)=L\}
\end{eqnarray*}
Since $\|{\bf x}-{\bf y}\|\geq \alpha(L)$ for ${\bf x}, {\bf y}\in L$ with ${\bf x}\neq {\bf y}$, we observe that $\{S_{\alpha(L)/2}({\bf a})\}_{{\bf a}\in K(L)}$ is a family of non-overlapping spheres touching the common sphere $S_{\alpha(L)/2}({\bf o})$. Therefore $|K(L)|$ is less than or equal to the maximum possible {\it kissing number} $k(d)$. It is known that $k(3)=12$,\footnote{This fact was conjectured correctly by Newton in a famous controversy between him and David Gregory (1694), and was proved by Sch\"{u}tte and van der Waerden in 1953.} and hence $|K(L)|\leq 12$ for a 3-dimensional lattice $L$. 
For example, 
\begin{eqnarray*}
&&\alpha(\mathbb{Z}^3)=1,~~K(\mathbb{Z}^3)=\{\pm{}^t(1,0,0),\pm{}^t(0,1,0),\pm{}^t(0,0,1) \},~~|K(\mathbb{Z}^3)|=6,\\
&&\alpha(L_{\mathcal{DT}})=\sqrt{3},~~K(L_{\mathcal{DT}})={\bf E}_{\mathcal{D}},~~|K(L_{\mathcal{DT}})|=8\\
&&\alpha(L_{\mathcal{D}})=\sqrt{2},~~K(L_{\mathcal{D}})={\bf E}_{\mathcal{DT}}~~|K(L_{\mathcal{D}})|=12.
\end{eqnarray*}
What we should notice is that $G(\mathbb{Z}^3)=G(L_{\mathcal{DT}})=G(L_{\mathcal{D}})$. Indeed, these groups coincide with the symmetry group ${\rm Iso}(Q)$ of the cube $Q$ introduced in the previous section, which acts transitively on the set of vertices not only of $Q$, but also of {\it cuboctahedron} and of the octahedron depicted in Fig.~\!\ref{octa3}.\footnote{The  cuboctahedron (also called the {\it heptaparallelohedron} or {\it dymaxion}) is one of thirteen {\it Archimedean solids}. See \cite{hisa} for interesting matters concerning cuboctahedron and the involvement of Archimedes, Pappus, Da Vinci, Kepler, etc. Perhaps Plato knew the cuboctahedron because Heron of Alexandria testifies in his {\it Definitiones} that Plato knew of a solid made of eight triangles and six squares.} Actually, ${\rm Iso}(Q)$ is identified with the full octahedral group ${\rm O_h}$, which acts irreducibly on $\mathbb{R}^3$.

\begin{figure}[htbp]
\vspace{-0.7cm}
\begin{center}
\hspace{0.7cm}
\includegraphics[width=.46\linewidth]{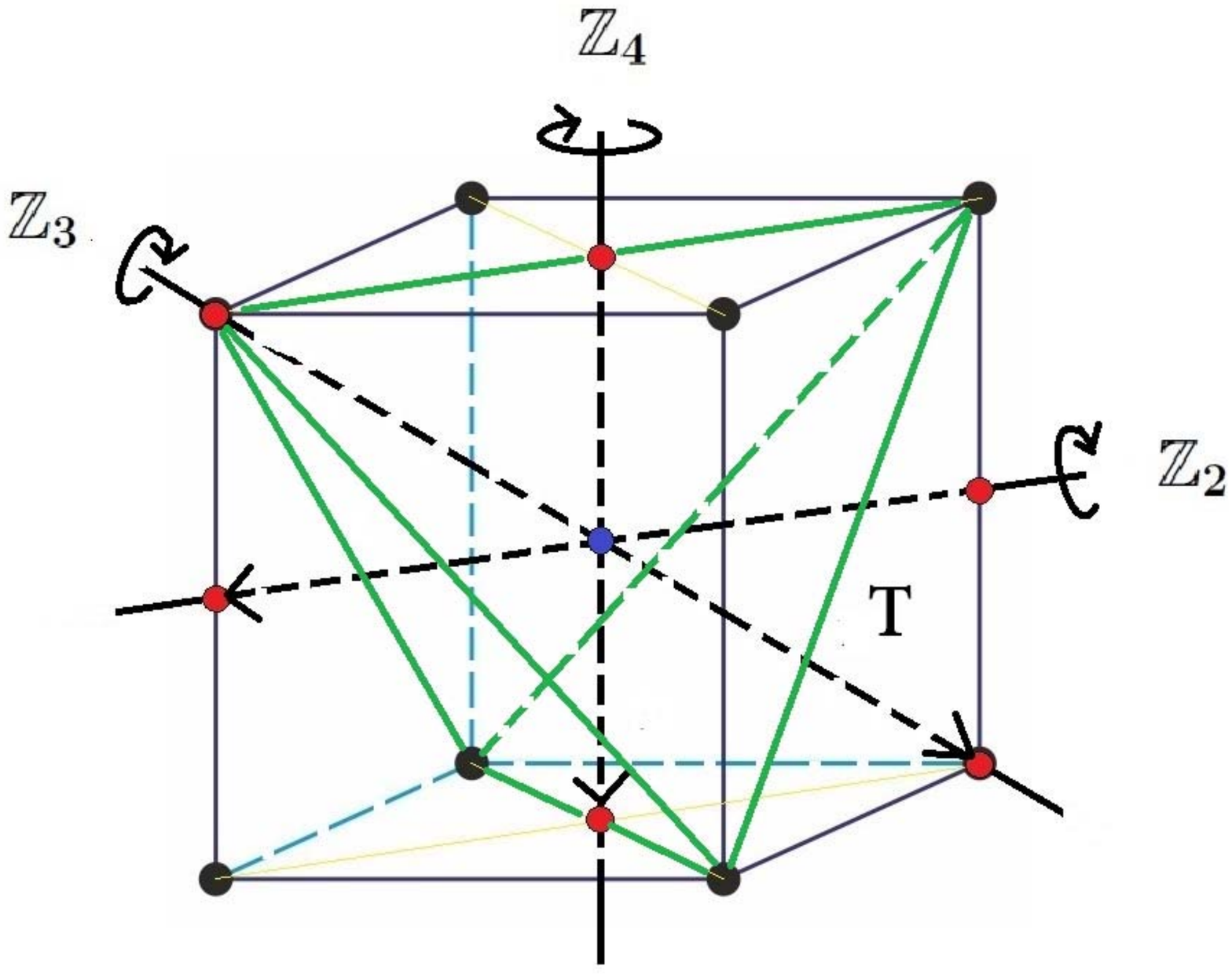}
\end{center}
\vspace{-0.8cm}
\caption{Octahedral rotations}\label{rotation}
\end{figure}

We shall mention a fact of importance about ${\rm Iso}(Q)$ which will be used later: 
The axis of a rotation in ${\rm Iso}(Q)$ is one of the thirteen lines passing through the origin and vertices of these three polyhedra. In other words, the family of lines 
$\mathbb{R}{\bf a}$
$({\bf a}\in K(L_{\mathcal{DT}})\cup K(L_{\mathcal{D}})\cup K(\mathbb{Z}^3))$ 
is exactly the family of axes of rotations in ${\rm Iso}(Q)$. In addition, there are three kinds of axes; that is, four axes passing through vertices ({\it the first kind}), six axes passing through the middle points of edges ({\it the second kind}), and three axes passing though the centroids of faces ({\it the third kind}). Representatives of these lines are depicted in Fig.~\!\ref{rotation}. Among the thirteen axes, four axes of the $1^{\rm st}$ kind and three axes of the $3^{\rm rd}$ kind are axes of rotations fixing the regular tetrahedron $T$ inscribed in the cube $Q$.

Now, keeping in mind the three examples of lattices $\mathbb{Z}^3$, $L_{\mathcal{DT}}$, $L_{\mathcal{D}}$, we make the following definition.

\begin{definition}{\rm 
A lattice $L$ in $\mathbb{R}^d$ is said to be {\it orthogonally symmetric} if

\smallskip
(i)~ 
$K(L)$ generates $L$, 

\smallskip
(ii)~ $G(L)$ acts transitively on $K(L)$, and

\smallskip

(iii)~ the $G(L)$-action on $\mathbb{R}^d$ is irreducible. 

}
\end{definition}
Notice that $G(L)=\{g\in {\rm O}(d)|~\!g(K(L))=K(L)\}$ because of Condition (i), and that $G(L)$ contains the {\it central reflection} $\sigma_0:(x,y,z)\mapsto (-x,-y,-z)$.

As observed above, the three lattices $\mathbb{Z}^3$, $L_{\mathcal{DT}}$, and $L_{\mathcal{D}}$ are orthogonally symmetric.

Typical examples of orthogonally symmetric lattices of general dimension are {\it irreducible root lattices} whose properties actually motivate the definition above. For the convenience of the reader, let us recall the definition of root lattices.  

For an {\it even lattice} $L$, i.e., $\|{\bf x}\|^2\in 2\mathbb{Z}$ for all ${\bf x}\in L$, we let
$R(L)=\{{\bf x}\in L|~\!\|{\bf x}\|^2=2\}$. An element ${\bf x}\in R(L)$ is called a {\it root}. Clearly $\alpha(L)=\sqrt{2}$ and $R(L)=K(L)$.

\begin{definition} {\rm 
 An even lattice $L$ is called a {\it root lattice} if the root system $R(L)$ generates $L$.
 A root lattice $L$ is said to be {\it irreducible} if $L$ is not a direct sum of two non-trivial lattices. 

} 
\end{definition}

Orthogonal symmetricity of irreducible root lattices is derived from the fact that the {\it Weyl group}, a subgroup of $G(L)$ generated by reflections through hyperplanes associated to the roots, acts transitively on $R(L)=K(L)$, and acts irreducibly on $\mathbb{R}^d$. Representatives of irreducible root lattices are $A_d$ and $D_d$ ($A_3=D_3$) in the usual notations for root systems.\footnote{Irreducible root lattices are classified into two infinite families of classical root lattices $A_d ~(d\geq 1), D_d ~(d\geq 4)$ and the three {\it exceptional} lattices $E_6, E_7, E_8$ (cf. \!\cite{ebe}). Another remarkable example of an orthogonally symmetric lattice is the {\it Leech lattice} in $\mathbb{R}^{24}$, symbolically denoted by $\Lambda_{24}$, discovered by John Leech in 1967 (cf. \cite{mil}), for which we have $\alpha(\Lambda_{24})=2$ and $|K(\Lambda_{24})|=196560$.} Here $A_d$ is the lattice in the orthogonal complement $(\underbrace{1,\ldots,1}_{d+1})^{\perp}$ in $\mathbb{R}^{d+1}$ defined by
$$
A_d=\{(x_1,\ldots, x_{d+1})\in \mathbb{Z}^{d+1}|~\!x_1+\cdots+x_{d+1}=0\}.
$$
The root system $R(A_4)$ consists of vectors such that all but two coordinates equal to $0$, one coordinate equal to $1$, and one equal to $–1$. As for $D_d$, it  is defined as
$$
D_d=\{(x_1,\ldots,x_d)\in \mathbb{Z}^d|~\!x_1+\cdots+x_d~\text{even}\}.
$$ 
The root system $R(D_d)$ consists of all integer vectors in $\mathbb{R}^d$ of length $\sqrt{2}$.

Note that the root lattice $A_3~(=D_3)$ coincides with $L_{\mathcal{D}}$.\footnote{The root lattice $A_d$ is the period lattice of the {\it $d$-dimensional diamond} that is defined as the standard realization of the maximal abelian covering graph of the dipole graph with $d+1$ parallel edges (\cite{kotani3}).}

\medskip

The following is the statement referred to in the previous section.

\begin{proposition}\label{prop:2}
The three lattices $\mathbb{Z}^3$, $L_{\mathcal{DT}}$ and $L_{\mathcal{D}}$ (up to similarity) are the only examples of orthogonally symmetric lattices of $3$-dimension.

\end{proposition}

\noindent{\it Proof} ~
Let $L$ be an orthogonally symmetric lattice. It suffices to show that $K(L)$ is one of $K(\mathbb{Z}^3)$, $K(L_{\mathcal{DT}})$. $K(L_{\mathcal{D}})$ (up to similarity).

The classification of finite subgroups of ${\rm O}(3)$ containing the central reflection tells us that the possibilities for $G(L)$ are exhausted by $\mathbb{Z}_k\times \mathbb{Z}_2$, $\mathcal{D}_k\times \mathbb{Z}_2$, 
$\mathcal{S}_4\times \mathbb{Z}_2$, $\mathcal{A}_4\times \mathbb{Z}_2$, $\mathcal{A}_5\times \mathbb{Z}_2$ (cf. \cite{New}), where, given a line $\ell$ passing through the origin, 
\begin{eqnarray*}
&&\mathbb{Z}_k=\text{cyclic group generated by the rotation $\rho$ of angle $2\pi/k$ with the axis $\ell$}\\
&&\mathcal{D}_k=\text{dihedral group generated by $\rho$ and a reflection in a plane containing $\ell$}, 
\end{eqnarray*}
and $\mathcal{A}_k$ is the alternating group of degree $k$ which, in the case $k=4$, acts on a regular tetrahedron by rotations. The direct factor $\mathbb{Z}_2$ is the group generated by the central reflection. But $\mathbb{Z}_k\times \mathbb{Z}_2$ and $\mathcal{D}_k\times \mathbb{Z}_2$ 
are obviously excluded. Furthermore, $\mathcal{A}_5\times \mathbb{Z}_2$, which is isomorphic to the {\it full icosahedral group}, is also excluded in view of the crystallographic restriction (cf. \cite{engel}).\footnote{If $G(L)$ contains an element of order $n$, then $\varphi(n)\leq d$, where $\varphi$ is the {\it Euler function}; in particular, $n=1,2,3,4,$ or $6$ in the 3-dimensional case.}

\begin{figure}[htbp]
\vspace{-0.6cm}
\begin{center}
\includegraphics[width=.54\linewidth]{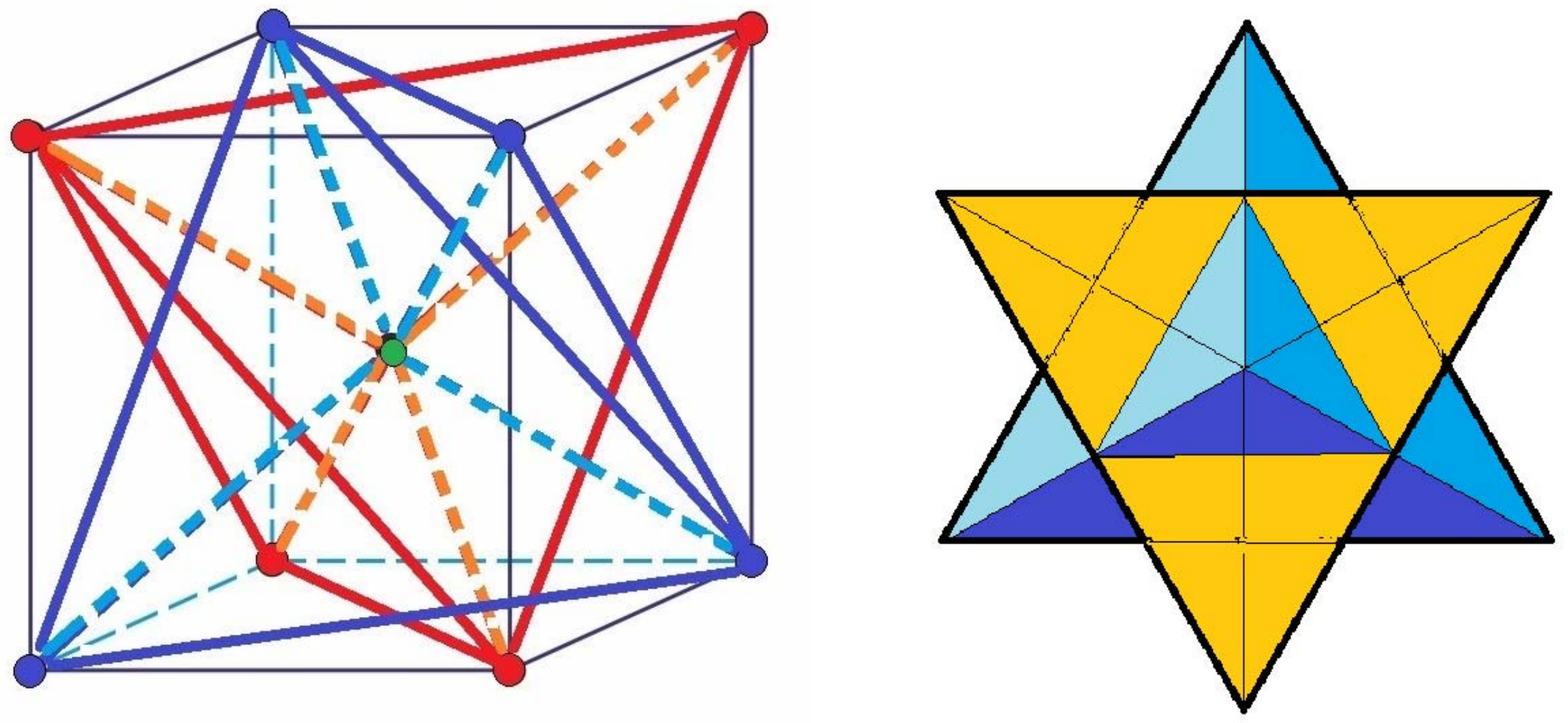}
\end{center}
\vspace{-2cm}
\caption{Star tetrahedron}\label{star}
\end{figure}

Take the cube $Q$ mentioned avove. 
If $G(L)\cong \mathcal{S}_4\times \mathbb{Z}_2$, then $G(L)$ is the full octahedral group, and is identified with ${\rm Iso}(Q)$. If $G(L)\cong \mathcal{A}_4\times \mathbb{Z}_2$, then $G(L)$ acts, in a natural manner, on a {\it star tetrahedron}, the solid composed of a tetrahedrons and its transformation by $\sigma_0$, and hence is assumed to acts on $Q$ as well (Fig.~\!\ref{star}). Thus in any case, $G(L)$ is identified with a subgroup of ${\rm Iso}(Q)$. Furthermore, even when $G(L)\cong \mathcal{A}_4\times \mathbb{Z}_2$,
the action of $G(L)$
on $Q$ 
is vertex-, edge-, and face-transitive, from which it follows that $G(L)$ also acts transitively on each set of vertices of the octahedron, cube and cuboctahedron.

Take an element ${\bf a}\in K(L)$, and consider the isotropy group $G_{\bf a}=\{g\in G(L)|~\!g{\bf a}={\bf a}\}$, which is isomorphic to $\mathbb{Z}_k$ or $\mathcal{D}_k$ acting on the plane which passes through the origin and is orthogonal to ${\bf a}$. 

We now have
$$
|K(L)|=2,4,6,8,10,\text{or}~12,
$$
Here the two cases $|K(L)|=2$ or $4$ are excluded for an obvious reason. The case $|K(L)|=10$ is also excluded since, otherwise, $|G|~(=|G_{\bf a}|\cdot |K(L)|)$ is divisible by $5$; this does not happen again in view of the crystallographic restriction.
Therefore $|K(L)|=6,8,12$. What remains to do is to determine the configuration of vectors in $K(L)$ in each case. This will be done by brute-force, rather than a systematic method.

In the case that $|K(L)|=6$, using $|\mathcal{A}_4\times \mathbb{Z}_2|=24$ and $|\mathcal{S}_4\times \mathbb{Z}_2|=48$, 
we have $|G_{\bf a}|=4$ or $8$, so the possible structures of $G_{\bf a}$ are described as
$$
G_{\bf a}\cong\begin{cases}
\mathbb{Z}_4~\text{or}~\mathcal{D}_2 & (|G_{\bf a}|=4;~ \text{the case $G(L)\cong\mathcal{A}_4\times \mathbb{Z}_2$})\\
\mathbb{Z}_8~\text{or}~\mathcal{D}_4 & (|G_{\bf a}|=8;~\text{the case $G(L)\cong\mathcal{S}_4\times \mathbb{Z}_2$})
\end{cases}.
$$ 
Thus in any case, $G_{\bf a}$ contains a rotation whose axis is $\mathbb{R}{\bf a}$. This axis must be one of the 13 lines mentioned before since $G_{\bf a}\subset {\rm Iso}(Q)$. Thus a non-zero scalar multiple of {\bf a}, say $c{\bf a}$, belongs to one of $K(\mathbb{Z}^3)$, $K(L_{\mathcal{DT}})$, $K(L_{\mathcal{D}})$, and $K(L)~(=G(L){\bf a}$) coincides with one of $c^{-1}K(\mathbb{Z}^3)$, $c^{-1}K(L_{\mathcal{DT}})$, $c^{-1}K(L_{\mathcal{D}})$, where we use the fact that $G(L)$ acts transitively on the set of vertices of the three polyhedron mentioned above. Since $|K(L)|=6$, we conclude that $K(L)=c^{-1}K(\mathbb{Z}^3)$. 

The same argument applies to the cases $|K(L)|=8$ because
$$
G_{\bf a}\cong\begin{cases}
\mathbb{Z}_3 & (|G_{\bf a}|=3;~ \text{the case $G(L)\cong\mathcal{A}_4\times \mathbb{Z}_2$})\\
\mathbb{Z}_6~\text{or}~\mathcal{D}_3 & (|G_{\bf a}|=6;~\text{the case $G(L)\cong\mathcal{S}_4\times \mathbb{Z}_2$})
\end{cases},
$$ 
and leads to the conclusion that there is a non-zero constant $c$ such that $K(L)=c^{-1}K(L_{\mathcal{DT}})$.
 
\begin{figure}[htbp]
\vspace{-0.7cm}
\begin{center}
\includegraphics[width=.85\linewidth]{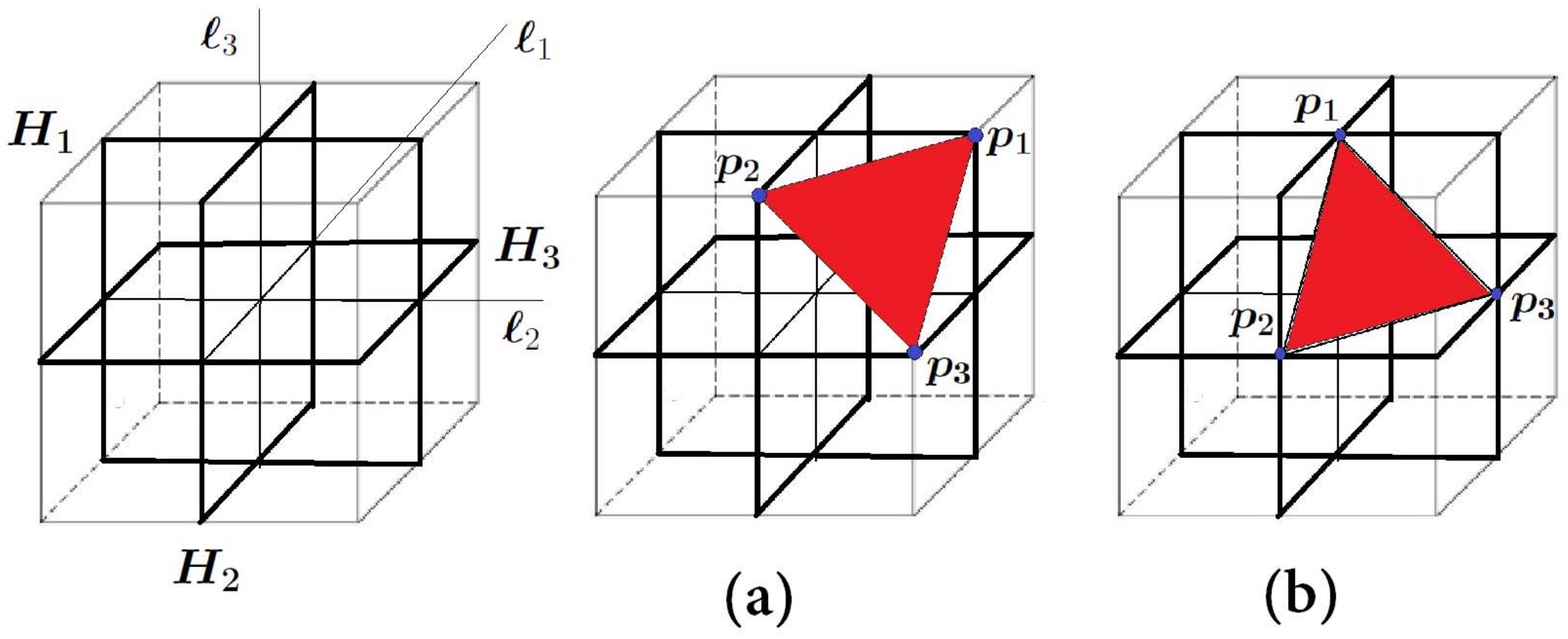}
\end{center}
\vspace{-4cm}
\caption{}\label{cube3}
\end{figure}

In the case that $|K(L)|=12$, we have $|G_{\bf a}|=2$ or $4$, so 
$$
G_{\bf a}\cong\begin{cases}
\mathbb{Z}_2~\text{or}~\mathcal{D}_1 & (|G_{\bf a}|=2;~ \text{the case $G(L)\cong\mathcal{A}_4\times \mathbb{Z}_2$})\\
\mathbb{Z}_4~\text{or}~\mathcal{D}_2 & (|G_{\bf a}|=4;~\text{the case $G(L)\cong\mathcal{S}_4\times \mathbb{Z}_2$})
\end{cases}.
$$ 
The possible occurrence of $G_{\bf a}=\mathcal{D}_1$ may cause a difficulty. If this really occurs, then 
$G_{\bf a}$ is generated by a reflection, say $\mu$,  
through a certain plane $H_{\bf a}$ containing ${\bf a}$,
which is expressed as $\rho\sigma_0$ with a rotation $\rho$ belonging to the symmetry group for the tetrahedron $T$
 since $G(L)\cong \mathcal{A}_4\times \mathbb{Z}_2$. From $\rho^2=\mu^2=1$, it follows that $\rho~(=\mu\sigma_0$) 
is a rotation of angle $\pi$
whose axis $\ell_{\bf a}$ (necessarily of the $3^{\rm rd}$ kind since rotations of the $1^{\rm st}$ kind are of order 3) is orthogonal to $H_{\bf a}$. Letting ${\bf a}$ run over $K(L)$ and ignoring the duplicate, we have three axes $\ell_1, \ell_2, \ell_3$ of 
the $3^{\rm rd}$ kind and also three planes $H_1, H_2, H_3$ orthogonal to each other (see the left of Fig.~\!\ref{cube3}), each of which contains four vectors in $K(L)$. We then take a look at the four rotations whose axes are of the $1^{\rm st}$ kind. They belong to $G(L)$ since $G(L)$ has the subgroup corresponding to $\mathcal{A}_4$, and each of them yields a cyclic permutation of $H_1, H_2, H_3$. Moreover, each orbit of such rotations acting on the set of end points of vectors in $K(L)$ form an equilateral triangle whose vertices $p_1,p_2,p_3$ are, we may assume, in the planes $H_1, H_2, H_3$, respectively. This imposes a restriction on the location of the triangle; indeed, there are just two possible cases of its location as depicted in Fig.~\!\ref{cube3} (a) and (b).\footnote{To make the argument simple, we assume that $p_i$'s are on the faces of the cube $Q$.} We perform this procedure for each axis. In case (b) (i.e., the case that $p_i$ lies in the intersection of two planes), we obtain the octahedron, and hence $|K(L)|=6$, contradicting to the assumption that $|K(L)|=12$. In case (a), we see with a little thought that $K(L)$ yields the cuboctahedron. This implies that $K(L)$ is similar to $K(L_{\mathcal{D}})$.\footnote{What was actually proved is that the case $G(L)\cong \mathcal{A}_4\times \mathbb{Z}_2$ does not occur.}

The proof of Proposition \ref{prop:2} is thus completed.

\begin{remark}\label{rem}{\rm 
(i)~ 
The root lattice $A_3~(=D_3=L_{\mathcal{D}})$ has something to do with one of the densest packings of equally sized spheres (called the {\it face-centered cubic packing}; Fig.~\!\ref{hexa}, where the nodes (vertices) are located at center of equal spheres, and line segments (edges) indicate that two spheres located at their end points touch each other).\footnote{It was Kepler who speculated that the densest packing is attained by the face-centered cubic packing ({\it Strena Seu de Nive Sexangula}, 1611). His conjecture, considered the earliest statement on crystal structures, was confirmed by T. Hales in 1998 with the aid of a computer (prior to this,  Gauss proved that the Kepler conjecture is true if restricted to {\it lattice packings} ({\it Werke}, III, 481--485). 
For $d$ with $4\leq d\leq 8$ and $d=24$, it is known that the densest lattice packing in $\mathbb{R}^d$ is attained by the irreducible root lattices $D_4$, $D_5$, $E_6$, $E_7$, $E_8$, and the Leech lattice, respectively (cf.~\!\cite{Con}), which are also lattices providing the maximal kissing numbers in each dimension.}

\begin{figure}[htbp]
\vspace{-0.7cm}
\begin{center}
\includegraphics[width=.7\linewidth]{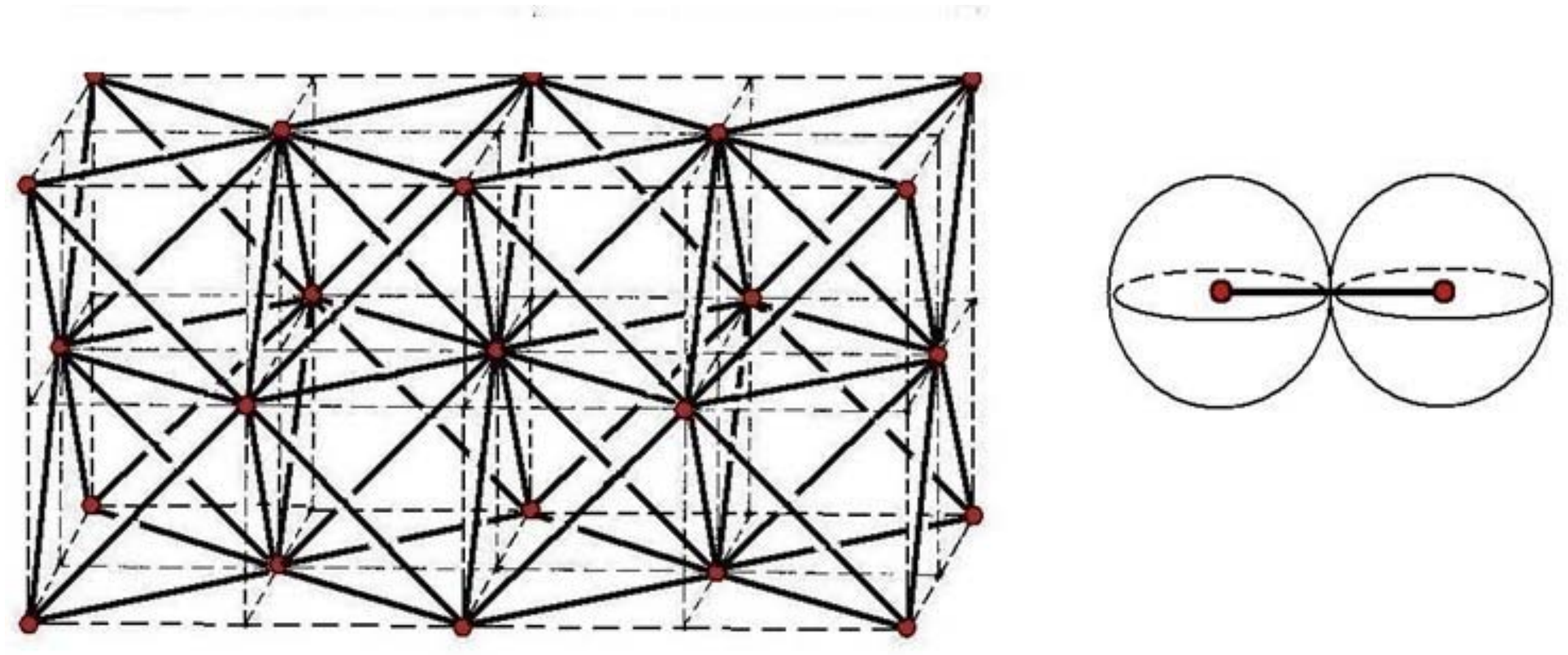}
\end{center}
\vspace{-2.9cm}
\caption{Face-centered cubic packing}\label{hexa}
\end{figure}

\medskip

(ii)~ Let $L$ be a lattice in $\mathbb{R}^d$ such that $K(L)$ generates $L$.
For two vectors ${\bf a}$ and ${\bf b}$ in $K(L)$ such that ${\bf b}\neq -{\bf a}$, the angle $\theta({\bf a},{\bf b})$ 
between ${\bf a}$ and ${\bf b}$ satisfies 
\begin{equation*}\label{eq:ineq}
\frac{\pi}{3} \leq \theta({\bf a},{\bf b})\leq \frac{2\pi}{3},
\end{equation*}
because otherwise $\|{\bf a}+{\bf b}\|<\alpha(L)$ or $\|{\bf a}-{\bf b}\|<\alpha(L)$.
 Using this, one can show that orthogonally symmetric lattices of $2$-dimension are either the square lattice or the regular triangular lattice (Fig.~\!\ref{regtri}).\footnote{As well-known, positive definite quadratic forms are associated with lattices. What we have observed tells us that, among all positive definite binary quadratic forms, the classical ones $x^2+y^2$ and $x^2\pm xy+y^2$ are distinguished from others by means of symmetry. These quadratic forms have played a significant role in the history of number theory,} 

\begin{figure}[htbp]
\vspace{-0.8cm}
\begin{center}
\includegraphics[width=.79\linewidth]{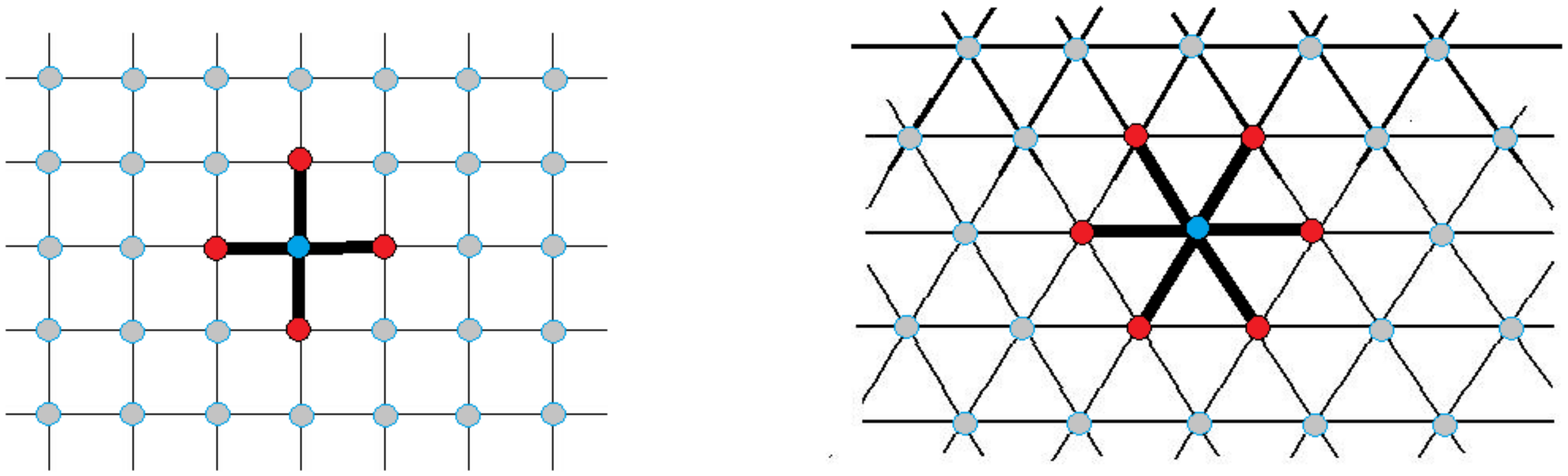}
\end{center}
\vspace{-4.2cm}
\caption{Orthogonally symmetric lattices of 2-dimension}\label{regtri}
\end{figure}


(iii)~  
If $L$ is orthogonally symmetric, then $K(L)$ is a {\it tight frame}, the notion appearing in wavelet analysis and crystal design (\cite{sunada1}). That is, 
there exists a positive constant $c$ such that
\begin{equation}\label{eq:11}
\sum_{{\bf a}\in K(L)}\langle{\bf a},{\bf x}\rangle{\bf a}=c{\bf x}\quad ({\bf x}\in \mathbb{R}^d).
\end{equation}
Explicitly
$$
c=\frac{1}{d}\alpha(L)^2|K(L)|.
$$
Indeed, letting ${\bf e}_i$ $(i=1,\ldots,d)$ be the standard basis, and applying (\ref{eq:11}) to ${\bf x}={\bf e}_i$, we find
$
c=\sum_{{\bf a}\in K(L)}\langle {\bf a},{\bf e}_i\rangle^2
$.
Since $\sum_{i=1}^d\langle {\bf a},{\bf e}_i\rangle^2=\|{\bf a}\|^2=\alpha(L)^2$, we obtain
$$
dc=\alpha(L)^2|K(L)|.
$$

}

\end{remark}

\end{document}